\documentclass[preprintnumbers,amsmath,amssymb,aps]{revtex4}

\usepackage{graphicx}

\let\le=\leqslant  \let\ge=\geqslant
\newcommand\PA{\partial}     \newcommand\DS{\displaystyle}
\newcommand\MR{\mathrm}
\newcommand\R{\mathbb{R}}     \newcommand\Z{\mathbb{Z}}

     \newcommand\Imag{\mbox{Im}}
\newcommand\vect[1]{{\bf #1}}     \newcommand\cri{cr}
\newcommand\dint[1]{\,\MR d#1}     \newcommand\Conj{\mbox{Conj}}
\newcommand\refsecc[1]{~\!\S\,\ref{#1}}
\newcommand\reffig[1]{figure~\ref{#1}}
\newcommand\refecua[1]{(\ref{#1})}
\newcommand\reftab[1]{table\ \ref{#1}}
\newcommand\figura[2]{\begin{figure}[t]\begin{center}
  \includegraphics{#1.pdf}\caption{#2.}
  \label{#1}\end{center}\end{figure}}
\newcommand\apartado[1]{\medbreak\noindent{\bf #1.}}
\providecommand\bnabla{\boldsymbol{\nabla}}
\providecommand\bcdot{\boldsymbol{\cdot}}

\begin{document}

\title{Hopf bifurcations to quasi-periodic solutions\\for the two-dimensional
plane Poiseuille flow}

\author{Pablo S. Casas}
 \homepage{http://www.ma1.upc.es/~casas}
\affiliation{Departamento de Matem\'atica Aplicada I, Universidad Polit\'ecnica
 de Catalu\~na, Diagonal, 647. 08028 Barcelona. Spain.}

\author{\`Angel Jorba}
 \homepage{http://www.maia.ub.es/~angel}
\affiliation{Departamento de Matem\'atica Aplicada y An\'alisis, Universidad de
Barcelona, Gran Via, 585. 08007 Barcelona. Spain.}

\date{\today}

\begin{abstract}
This paper studies various Hopf bifurcations in the two-dimensional plane Poiseuille
problem. For several values of the wavenumber $\alpha$, we obtain the branch of
periodic flows which are born at the Hopf bifurcation of the laminar flow. It is
known that, taking $\alpha\approx1$, the branch of periodic solutions has several
Hopf bifurcations to quasi-periodic orbits. For the first bifurcation, previous
calculations seem to indicate that the bifurcating quasi-periodic flows are stable
and go backwards with respect to the Reynolds number, $Re$. By improving the
precision of previous works we find that the bifurcating flows are unstable and go
forward with respect to $Re$. We have also analysed the second Hopf bifurcation of
periodic orbits for several $\alpha$, to find again quasi-periodic solutions with
increasing $Re$. In this case the bifurcated solutions are stable to superharmonic
disturbances for $Re$ up to another new Hopf bifurcation to a family of stable
$3$-tori. The proposed numerical scheme is based on a full numerical integration of
the Navier-Stokes equations, together with a division by 3 of their total dimension,
and the use of a pseudo-Newton method on suitable Poincar\'e sections. The most
intensive part of the computations has been performed in parallel. We believe that
this methodology can also be applied to similar problems.
\end{abstract}

\maketitle

\section{Introduction}

The theory of hydrodynamic stability is one of the main topics in fluid mechanics.
Poiseuille as well as Taylor--Couette flow are test problems where it is possible
the evaluation of different analytical and numerical methods, due essentially to the
simplicity of their geometry. The dynamics of plane Poiseuille flow departs from the
laminar flow. The stability of the laminar solution to infinitesimal disturbances
has been analysed linearly and gives rise to the Orr--Sommerfeld equation. This
equation has been studied by several authors as \citet{Thomas53}, \citet{Orszag71},
and \citet{Maslowe85} among others, and it is well understood. The critical Reynolds
number of the linear theory, $Re_{\cri}=5772.22$ for the wavenumber
$\alpha=1.02056$, has been obtained by this approach. However, as experiments of
\citefullauthor{Carlson82} \cite{Carlson82}, \citet{Nishioka85}, and
\citefullauthor{Alavyoon86} \cite{Alavyoon86} showed, transition to turbulence is
observed for Reynolds number $\approx1000$, what motivates that finite-amplitude
disturbances originate the transition. The understanding of the transition to
turbulence has been conjectured by \citet{Saffman83} to depend on intermediate
vortical states and turbulence takes place due to their three-dimensional
instability. In recent years, authors have also payed attention to subcritical
transition models based on transient optimal growth (see \citet{Schmid01}, for
instance). Examples of vortical states are periodic\footnote{Unless stated otherwise
  ``periodic'' or ``quasi-periodic'' refers to time in a fixed frame of reference.}
flows in time or space, among which can be mentioned: two-dimensional travelling
waves, secondary flows in two or three dimensions (for them the flow rate and the
pressure gradient are constants) and quasi-periodic solutions. \citet{Ehrenstein91}
discovered a new family of secondary bifurcation branches in dimension 3, which
contains only even spanwise Fourier modes and reduces the critical Reynolds number
(defined in terms of the averaged velocity across the channel) to
$Re_{Q_m}\approx1000$ as observed in experiments.

Two-dimensional disordered motion is associated with the large scales of some
turbulent flows, so there probably exist attractors for those two-dimensional flows.
Besides, two- and three-dimensional states can compete and coexist in the final flow
(cf. \citet{Jiménez87} and the references therein). In spite of the fact that
transition to turbulence is a three-dimensional phenomenon, there are many
properties of the two-dimensional flows observed in fully turbulent
three-dimensional flows such as wall sweeps, ejections, intermittency and bursting,
as \citet{Jiménez90} showed. The two-dimensional case has attracted the attention of
many authors but it is not completely understood as the problem of two-dimensional
transition to turbulence proves. Due to \citeauthor{Squire33}'s \cite{Squire33}
theorem, to every three-dimensional perturbation of the linearized Navier--Stokes
equations for a given $Re,\alpha$, it corresponds a two-dimensional one for some
$\tilde\alpha\ge\alpha$ and $\tilde{Re}\le Re$, so the critical $Re$ for the linear
theory must be attained by a two-dimensional flow. This result has been one of the
main reasons to firstly try to understand the two-dimensional case, apart from the
obvious easiness of computations compared to the three-dimensional situation. In
addition, some of the properties obtained from the two-dimensional case can also
provide new insight for three-dimensional flows.

In this work we intend to analyse the dynamics of an easily treatable problem
without domain complexities as is the case of the two-dimensional plane Poiseuille
flow. Different levels of bifurcation to respective vortical states are considered,
starting at the basic parabolic flow. From it, a family of travelling waves is born
subcritically (see \refsecc{estab_sol_per}) for $\alpha\approx1$. There are many
papers concerning this kind of waves: \citet{Soibelman91} gave an excellent review
about it. As a starting point for our computations we have also reproduced the
calculations to find the travelling waves for several values of $\alpha$.
\citet{Jiménez87,Jiménez90} and \citet{Soibelman91} obtained the next level of
bifurcation to quasi-periodic solutions. Employing full numerical simulation in
time, \citet{Jiménez87,Jiménez90} computed different attractor flows with a moderate
number of Chebyshev and Fourier modes. On the other hand, \citet{Soibelman91}
implemented an algebraic approach to capture stable and unstable quasi-periodic
flows, but the number of modes used were not enough to give good results and they
were not able to carry out the stability analysis. The method implemented in the
present work combines both: we solve a stationary problem to compute travelling
waves for an observer moving at an appropriate speed, whereas the quasi-periodic
flows are found by means of full numerical integration of the Navier--Stokes
equations. Through algebraic manipulations, we express the discretized
Navier--Stokes system only in terms of the stream component of the velocity. As a
consequence, the dimension of the system is divided by 3, reducing considerably the
computational effort. Using the numerical integrator, we have built a Poincar\'e
section of the flow, in order to apply a pseudo-Newton method for obtaining also
unstable quasi-periodic solutions. These unstable intermediate states of the flow
provide a highly useful insight into the transition process, as exemplified by
secondary bifurcations in shear flows (see \citet{Casas04} for instance). The
spatio-temporal symmetries of the channel allows the reduction of quasi-periodic
flows with two-frequencies to periodic flows in the appropriate Galilean reference.
The quasi-periodic solutions found in this work correspond to the first two Hopf
bifurcations of travelling waves for the case of constant pressure drop through the
channel, and the first Hopf bifurcation when the mass flux is held constant. The
property of behaving as time-periodic flows if we take a suitable Galilean
reference, simplifies enormously the search of this kind of solutions. For them, the
associated return time to the Poincar\'e section is roughly $10000$ time units at
the first Hopf bifurcation for constant pressure, what makes the temporal
integration very costly. The considered numerical procedure utilizes a parallel
algorithm to evaluate the different columns of a Jacobian matrix, needed in the
application of pseudo-Newton's method for the continuation of quasi-periodic
solutions. We find that on the analysed Hopf bifurcations for both constant pressure
and constant flux formulations, there exist quasi-periodic flows with increasing
$Re$ for some range of $\alpha$ and with decreasing $Re$ for some other $\alpha$:
the bifurcations are supercritical or subcritical respectively. On the first
bifurcation for constant pressure, we have traversed a curve of unstable
quasi-periodic solutions. On the remaining bifurcations, there are stable
quasi-periodic solutions to disturbances with the same wavenumber $\alpha$ and
likewise, for $Re$ sufficiently large, we have obtained unstable solutions.

Once we have situated the different studies concerning Poiseuille flow, in the next
section we pose the concrete terms that define the plane Poiseuille problem in two
dimensions, together with their equations for both cases of constant pressure and
flux. Next in\refsecc{numerical_app} we explain the main details of the numerical
methods. In\refsecc{solper} we review some results of the Orr--Sommerfeld equation
and obtain, for several values of $\alpha$, the bifurcating solutions of
time-periodic flows. From these we analyse in\refsecc{solcper} the bifurcating
branches to quasi-periodic solutions at the above-mentioned Hopf bifurcations.
Finally in\refsecc{conclusions} we point out some conclusions.

\section{Poiseuille flow}\label{Pois-flow}

We consider the flow of a viscous incompressible two-dimensional fluid, in a
channel between two parallel walls, governed by the Navier--Stokes equations
together with the incompressibility condition
\begin{equation}
\DS\frac{\PA\vect u}{\PA t} + (\vect u\bcdot\bnabla)\vect u =-\nabla p +
\frac{1}{Re}\Delta\vect u, \qquad\bnabla\bcdot \vect u = 0, \label{NS}
\end{equation} where $\vect
u=\vect u(x,y,t)=(u,v)(x,y,t)$ represents the two-dimensional velocity, $p=p(x,y,t)$
the pressure and $Re$ the Reynolds number. As boundary conditions we suppose no-slip
on the channel walls at $y=\pm 1$ and, at artificial boundaries in the stream
direction $x$, a period $L$, i.e.
\begin{equation}\left.\begin{matrix}
\hfill u(x,\pm1,t)=v(x,\pm1,t)=0\cr \hfill
(u,v,p')(x+L,y,t)=(u,v,p')(x,y,t)\cr
\end{matrix}\right\} \qquad x\in\R, \quad y\in[-1,1],
\quad t\ge0, \label{CF}\end{equation} being $p'=p+Gx$, for $G = G(t)$ the mean
pressure gradient on the channel length, $L$, in the streamwise direction. For the
system previously described there is a time-independent solution known as the basic
or laminar flow that has a parabolic profile, namely
\[u_b(y)=1-y^2, \quad v_b=0, \quad \nabla p_b=(-\frac{2}{Re},0).\]

Magnitudes in \refecua{NS}--\refecua{CF} are non-dimensional. We consider the two
typical formulations used to drive the fluid: fixing the total flux $Q$, or the mean
pressure gradient $G$, through the channel. For each of them we obtain a different
definition of $Re=hU_c / \nu$ namely, $Re_Q = 3Q / 4\nu$ and $Re_p = Gh^3 /
2\rho\nu^2$ respectively, where, in dimensional magnitudes, $h$ represents half of
the channel height, $U_c$ the velocity of the laminar flow in the centre of the
channel, and $\nu$ and $\rho$ the constant kinematic viscosity and density. For a
given laminar flow, i.e. letting $U_c$ fixed, both definitions of the Reynolds
number coincides with $Re=hU_c/\nu$. That is not the case for secondary flows,
defined as the ones having constant flux and mean pressure gradient through the
channel. If we consider such a flow $\vect u(x,y)$, expressed for each formulation
by means of respective Fourier series
\[\vect u^Q(x,y) = \sum_{k\in\Z}\hat{\vect
u}_k^Q(y)\MR{e}^{\MR{i}k\alpha x}, \qquad \vect u^p(x,y) = \sum_{k\in\Z}\hat{\vect
u}_k^p(y)\MR{e}^{\MR{i}k\alpha x},\] then, using the notation $[f]_a^b:=f(b)-f(a)$,
it is easy to check that (see for instance \citet{casas02})
\begin{equation} \frac{Re_p}{Re_Q} = -\frac{1}{4} \left[\frac{\PA{\hat u}_0^Q}{\PA
y}\right]_{-1}^1, \qquad \frac{Re_Q}{Re_p} = \frac{3}{4}\int_{-1}^1 u^p(x,y)\dint
y, \label{RelRey}\end{equation}
and the corresponding relationships between velocities and pressures
\begin{equation} \vect u^p(x,y) = \frac{Re_Q}{Re_p} \vect u^Q(x,y), \qquad p^p(x,y)
= \frac{Re_Q^2}{Re_p^2} p^Q(x,y). \label{RelVel}\end{equation}

We will employ later that periodic conditions at artificial boundaries in the stream
direction, yield a great simplification in the structure of the flow: quasi-periodic
solutions may be viewed as periodic flows, and periodic solutions as stationary
ones, if the observer moves at adequate speed $c$, in the stream direction. For this
reason we perform the change of variable $\tilde x=x-ct$, which (writing again $x$
instead of $\tilde x$) turns system \eqref{NS} into:
\begin{equation}\left\{
\begin{array}{r@{\hspace{0,5ex}}l}
\DS\frac{\PA u}{\PA t} + (u-c)\frac{\PA u}{\PA x} + v\frac{\PA u}{\PA y} &=
\DS-\frac{\PA p}{\PA x} + \frac{1}{Re}\left(\frac{\PA^2 u}{\PA x^2} +
\frac{\PA^2 u}{\PA y^2}\right) \\[3ex] \DS\frac{\PA v}{\PA t} + (u-c)\frac{\PA
v}{\PA x} + v\frac{\PA v}{\PA y} &= \DS-\frac{\PA p}{\PA y} +
\frac{1}{Re}\left(\frac{\PA^2 v}{\PA x^2} + \frac{\PA^2 v}{\PA y^2}\right) \\
[3ex] \DS\frac{\PA u}{\PA x} + \frac{\PA v}{\PA y} &=0,
\end{array}\right.\label{NSc}\end{equation} together with boundary conditions as
in \eqref{CF}. We can recover \eqref{NS} by simply taking $c=0$ in
\eqref{NSc}.

\section{Numerical approach}\label{numerical_app}

Let us now describe the numerical procedure. For system \eqref{NSc} we want to
follow the temporal evolution of an initial flow subjected to the incompressibility
condition, $\bnabla\bcdot\vect u=0$, and boundary conditions \eqref{CF}. To this end
we use a spectral method to approximate velocities $u,v$ and pressure deviation
$p'$, which from now on we consider non-dimensional quantities. We recall that
$p=p'-Gx$ and as it is easily obtained (see for example \citet{casas02})
\begin{equation}G = -\frac{1}{2Re_Q}\left[\frac{\PA\hat u_0}{\PA y}\right]_{-1}^1
\quad\mbox{or}\quad G=\frac{2}{Re_p},\label{G}\end{equation} respectively for the
constant flux or pressure cases, so in the first one the mean pressure gradient
varies with time and it is constant for the second one.

\apartado{Spatial discretization} We use a standard Fourier-Galerkin,
Chebyshev-collocation approach (cf. \citet{Canuto88}) in order to discretize $x,y$
derivatives. In this way, we consider Fourier series (with $\alpha=2\pi/L$ the
parameter wavenumber):
\[(u,v,p')(x,y,t) = \sum_{k=-N}^N(\hat u_k,\hat v_k,\hat p_k)
(y,t)\MR{e}^{\MR{i}k\alpha x}, \qquad x\in\R, \quad y\in[-1,1], \quad t\ge0,\] which
substituted in \eqref{NSc} gives rise to a system of partial differential equations
for the Fourier coefficients $(\hat u_k,\hat v_k,\hat p_k)$,
\begin{equation}\left\{
\begin{array}{r@{\hspace{0,5ex}}l}\DS\frac{\PA\hat u_k}{\PA t} +
\widehat{\left[(u-c)\frac{\PA u}{\PA x} + v\frac{\PA u}{\PA y}\right]_k} &=
\DS-\MR{i}k\alpha\hat p_k + \frac{1}{Re} \left(-k^2\alpha^2\hat u_k +
\frac{\PA^2\hat u_k}{\PA y^2} \right) + \delta_{k0}G, \\[3ex] \DS\frac{\PA{\hat
v_k}}{\PA t} + \widehat{\left[(u-c)\frac{\PA v}{\PA x} + v\frac{\PA v}{\PA
y}\right]_k} &= \DS-\frac{\PA\hat p_k}{\PA y} + \frac{1}{Re}
\left(-k^2\alpha^2\hat v_k+\frac{\PA^2\hat v_k}{\PA y^2}\right), \\[3ex]
\DS\MR{i}k\alpha\hat u_k + \frac{\PA\hat v_k}{\PA y} &= 0, \end{array}\right.
\label{NSdisc} \end{equation} where $-N\le k\le N,\ \widehat{[\cdot]}_k$ stands
for the order $k\mbox{th}$ Fourier coefficient of $[\cdot],\ \delta_{00}=1$, and
$\delta_{k0}=0$ for $k\neq0$. Because $u,v,p'$ are supposed to be real functions, it
is enough to consider modes $\hat u_k,\hat v_k,\hat p_k$ for $k=0,\ldots,N$ in
\eqref{NSdisc}. The corresponding no slip boundary conditions in \eqref{CF} are
now written as
\begin{equation}(\hat u_k,\hat v_k)(\pm1,t)=0, \quad\mbox{for } t\ge0\mbox{ and }
k=0,\ldots,N.\label{CFD}\end{equation}

The previous system is imposed at two different sets of Chebyshev abscissas to avoid
indeterminacy, namely $y_m=\cos(\pi m/M)$ (velocities and momentum) for $m = 1,
\dots, M-1$, and $y_{m+1/2}=\cos(\pi (m+1/2)/M)$ (pressure and continuity) for $m =
0, \dots, M-1$.

\apartado{Reduced equations} To emphasize the linear character of some operations,
we now write system \eqref{NSdisc} as
\begin{subequations}\begin{eqnarray} &&\dot u_k= -\left[(u-c)\frac{\PA u}{\PA x} +
v\frac{\PA u}{\PA y}\right]_k -D_{xk}C_1^{-1}C_2p_k + \frac{1}{Re} (D_{xk}^2 +
C_1^{-1}D_y^2C_1)u_k + \delta_{k0}G, \quad \label{NSmat1} \\ &&\dot v_k =
-\left[(u-c)\frac{\PA v}{\PA x} + v\frac{\PA v}{\PA y}\right]_k -C_1^{-1}D_yC_2p_k
+ \frac{1}{Re} (D_{xk}^2 + C_1^{-1}D_y^2C_1)v_k, \label{NSmat2} \\
&&D_{xk}C_2^{-1}C_1u_k + C_2^{-1}D_yC_1v_k = 0,\hfill \label{NSmat3}
\end{eqnarray}\label{NSmat}\end{subequations}
for $k=0,\ldots,N$, where we have taken `$\,\widehat{\phantom{\mbox{$u$}}}\,$' out
of $\widehat{[\cdot]}_k, \hat u_k,\hat v_k,\hat p_k$ for convenience. In
\refecua{NSmat} we have represented vectors of values $u_k, v_k$ at the grid $y_m$
and $p_k$ at the grid $y_{m+1/2}$; $C_1, C_2$ are the corresponding matrices of
cosines transforms for grids $y_m$ and $y_{m+1/2}$, and $D_{xk}, D_y$ denote the
respective matrices of partial derivatives in $x, y$.

From \eqref{NSmat3} we obtain a matrix $T_k$ that carries out the transformation
$\bar v_k = T_k\bar u_k$ where $\bar u_k=(u_{k,1},\ldots,u_{k,M-2})^t$ and $\bar
v_k=(u_{k,M-1}, v_{k,1},\ldots,v_{k,M-1})^t$ for $k=1,\ldots,N$. For $k=0$, from the
continuity equation in \eqref{NSdisc}, we obtain $\PA v_0/\PA y=0$. Applying
boundary conditions, $v_0(\pm1)=0$, we get $v_0(y)=0$. This implies
$v_{0,1}=\cdots=v_{0,M-1}=0$.

For $k=1, \dots, N$ we introduce the notation
\begin{eqnarray*}
&& \begin{aligned}
U_k &= -\left[(u-c)\frac{\PA u}{\PA x} + v\frac{\PA u}{\PA y}\right]_k +
\frac{1}{Re} (D_{xk}^2 + C_1^{-1}D_y^2C_1)u_k + \delta_{k0}G, \\
V_k &= -\left[(u-c)\frac{\PA v}{\PA x} + v\frac{\PA v}{\PA y}\right]_k +
\frac{1}{Re} (D_{xk}^2 + C_1^{-1}D_y^2C_1)v_k,
\end{aligned} \\
&& \begin{aligned}
\bar{U}_k &= (U_k)_{\{1,\dots,M-2\}}, \\
\bar{V}_k &= \begin{pmatrix}(U_k)_{\{M-1\}}\\ V_k\end{pmatrix},
\end{aligned}
\qquad\qquad
\begin{aligned}
\bar{Q}_k &= (D_{xk} C_1^{-1} C_2)_{\{1,\dots,M-2\}}, \\
Q_k &= \begin{pmatrix}(D_{xk} C_1^{-1} C_2)_{\{M-1\}}\\ C_1^{-1} D_y
C_2\end{pmatrix},
\end{aligned}
\end{eqnarray*}
where $A_{\{i_1, \dots, i_n\}}$ stands for rows $i_1, \dots, i_n$ of matrix $A$.
Equations \eqref{NSmat1} and \eqref{NSmat2} can be now expressed as
\[\begin{cases}
\dot{\bar u}_k = \bar U_k - \bar Q_k p_k, &\\ \dot{\bar v}_k = \bar V_k - Q_k p_k.
&\\ \end{cases}\]
The matrix $Q_k$ turns out to be an $M \times M$ invertible matrix. Consequently,
from the second equation we obtain $p_k = Q_k^{-1}( \bar V_k - \dot{\bar v}_k )$,
which substituted into the first one yields
\[\dot{\bar u}_k = \bar U_k - \bar Q_k Q_k^{-1}( \bar V_k - \dot{\bar v}_k ) =
\bar U_k - \bar Q_k Q_k^{-1}( \bar V_k - T_k\dot{\bar u}_k ).\]
Finally letting $P_k = \bar Q_k Q_k^{-1}$, we can also invert $I - P_k T_k$, and
thus we may solve for $\dot{\bar u}_k$
\begin{equation}\begin{cases} \dot u_0 = U_0, \\ \dot{\bar u}_k = (I - P_k
T_k)^{-1}(\bar U_k - P_k\bar V_k), & k=1, \dots, N, \end{cases}
\label{NSred}\end{equation} where $I$ is the identity matrix of dimension ${M-2}$
and we have extended the definition of $U_k$ for $k=0$. Bearing in mind the
substitution $\bar v_k = T_k\bar u_k$, we observe that system \eqref{NSred} does not
depend on $\bar v_k$ nor $p_k$: it only depends on $u_0$ and $\bar u_k$ for
$k=1,\ldots,N$. In addition, due to the elimination of pressure in \eqref{NSred}, we
avoid the indeterminacy caused by an additive constant. However this indeterminacy
has no effect upon the pressure gradient. Likewise this formulation saves the
problems in the imposition of consistent initial conditions with the
incompressibility. At the same time the stability analysis is simplified from
\eqref{NSred}.

\apartado{Temporal evolution} Once removed $v$ and $p$ from \eqref{NSmat}, in
\eqref{NSred} it just remains to discretize temporal derivatives. We can express
\eqref{NSred} as
\begin{equation}\dot{\bar u}_k = {\cal L}_k({\bar u}_k) + {\cal N}_k(\bar
u_0,\ldots,\bar u_N), \quad k=0, \dots, N, \label{NSredcomp}\end{equation} where
$\bar u_0 = u_0$ and ${\cal L}_k$, ${\cal N}_k$ corresponds respectively to linear
and nonlinear terms in $\bar u_0, \dots, \bar u_N$ on the right hand side of
\eqref{NSred}. We adopt a usual scheme for advection-diffusion problems: letting
$\bar u_k^n$ be $\bar u_k$ at the time instant $n\Delta t$ for some fix time step
$\Delta t$, we approximate ${\cal N}_k^j = {\cal N}_k(\bar u_0^j, \dots, \bar
u_N^j)$ by an explicit method (Adams--Bashforth) and ${\cal L}_k({\bar u}_k^n)$ by
an implicit one (Crank--Nicolson), so that \eqref{NSredcomp} yields
\begin{equation}\bar u_k^{n+1} - \frac{\Delta t}{2} {\cal L}_k(\bar u_k^{n+1}) =
\bar u_k^n + \frac{\Delta t}{2} \left[ {\cal L}_k(\bar u_k^n) + 3{\cal N}_k^n -
  {\cal N}_k^{n-1} \right].\label{IN}\end{equation} For the kind of solutions
treated in this work and moderate values of $Re \lesssim 10000$, we have verified
local errors originated in \eqref{IN} from the time discretization. For that
purpose, we approximate temporal derivatives by central finite differences and then
improve precision by means of extrapolations. In all tested cases we have found
errors $O((\Delta t)^2)$, which is in agreement with the discretization errors in
\eqref{IN}. For some flows considered in \refsecc{solcper} it has been necessary to
reduce $\Delta t$ to avoid overflows in $u(t)$.

We apply system \eqref{IN} to the two formulations described in \refsecc{Pois-flow},
namely, constant flux and constant mean pressure gradient. The imposition of
constant flux $Q = 4/3$ (a linear condition) allows us to reduce by one the number
of unknowns in $\bar u_0$. Therefore the number of equations is also reduced by one.
This condition is related to the formula derived for $G$ in \eqref{G}, which depends
linearly on $\bar u_0$ and thus it is included in ${\cal L}_0$. On the other hand,
in the constant pressure case, the value of $G$ is held constant and so it is a
nonlinear term. Taking into account that in \eqref{NSred}, $u_0$ has only real
components but for $k=1,\ldots,N,\ \bar u_k$ it is a complex vector, we conclude
that the block for $k=0$ has dimension $M-2$ or $M-1$ respectively for $Re_Q$ and
$Re_p$ formulations, and dimension $M-2$ for the $2N$ remaining real blocks. In
summary, each time step, \eqref{IN} implies the solution of a block diagonal linear
system of total real dimension $(2N+1)(M-2)$ in the constant flux case and
$(2N+1)(M-2)+1$ in the constant pressure one. That means a rough division by $3$ in
the dimension of the whole system \eqref{NSdisc}. In what follows we denote a
solution flow at time $t$ as $U(t) = (\bar u_0$, $\dots$,$\bar u_N)(t) \in \R^K$ for
$K = (2N+1)(M-2)+1$ or $K = (2N+1)(M-2)$, according to the two above-mentioned
cases.

\figura{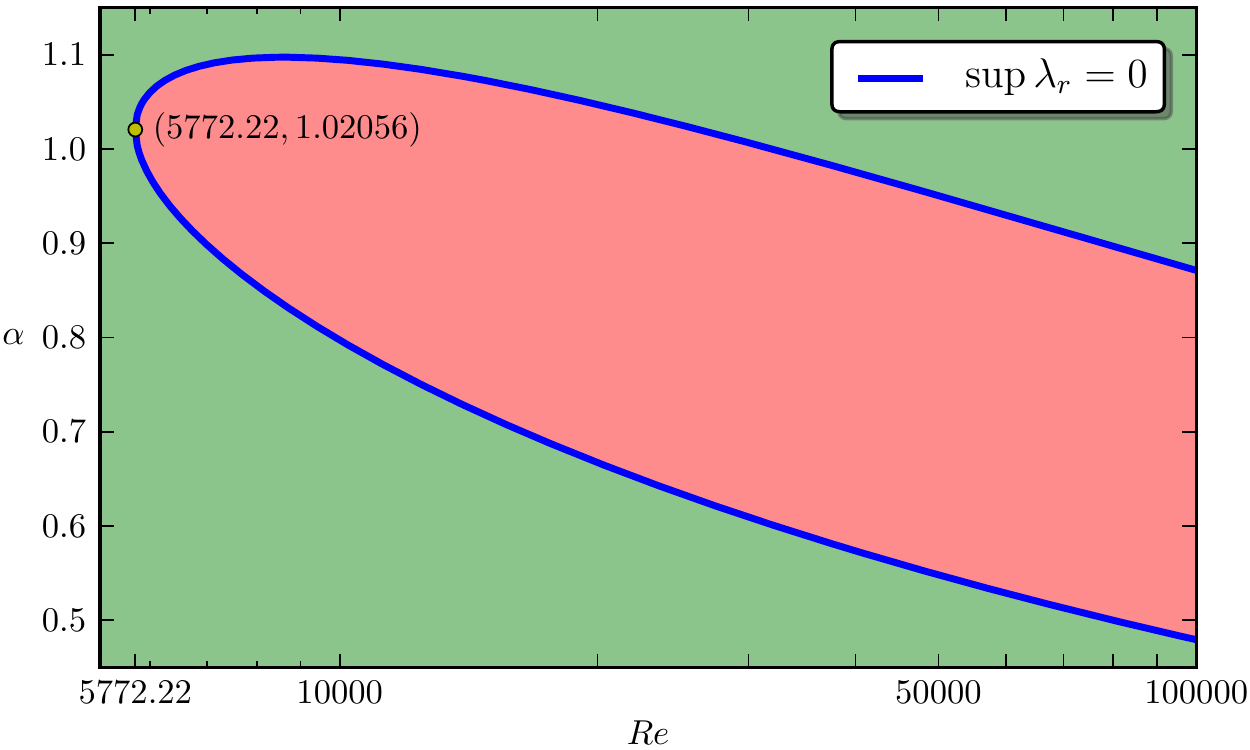}{Neutral stability curve (in blue) for the laminar solution,
using $n\lesssim1000$ discretization points. For each pair $(Re,\alpha)$ in this
curve, the most unstable eigenvalue $\lambda$ is purely imaginary. The curve splits
the $Re$-$\alpha$ plane in two stability regions as shown in the graph: the green
one is stable and the red one unstable}

\section{Periodic solutions}\label{solper}

\subsection{The Orr--Sommerfeld equation}

Before applying the previously described numerical scheme, we make some
considerations about the linearized stability of the laminar flow and time-periodic
solutions. We start from the linearization of the vorticity equation around the
basic flow, which is known as the Orr--Sommerfeld equation
\begin{equation}(u_b - \frac{\lambda\MR i}{\alpha})(\phi''-\alpha^2\phi)-u_b''\phi =
\frac{1}{\MR{i}\alpha Re}(\phi^{(4)} - 2\alpha^2\phi'' + \alpha^4\phi).
\label{OS}\end{equation} It is a fourth order ordinary differential equation on
$\phi = \phi(y)$ as eigenfunction, with $\lambda$ as eigenvalue, and boundary
conditions $\phi(\pm1)=\phi'(\pm1)=0$. For each $Re$ and $\alpha$, \eqref{OS}
represents an eigenvalue problem on $\lambda$ and $\phi$. In this way if $\lambda =
\lambda_r+\MR i\lambda_i$ is a complex eigenvalue with $\lambda_r>0$, then the
laminar flow is unstable to infinitesimal disturbances according to the linear
theory.

We have employed finite differences to approximate $\phi(y)$ and its derivatives in
an uniform mesh $\bar y_m=2m/(n+1)-1\in[-1,1]$ for $m=0,\ldots,n+1$ and $n$ a
sufficiently large positive integer. After substituting $\phi(\bar
y_m),\ m=0,\ldots,n+1$ and the approximation to its derivatives in \eqref{OS}, we
obtain an eigenvalue problem of finite dimension: $A\phi = cB\phi$, for $A,B$
matrices depending only on $Re,\alpha$, and $n$: $A$ is pentadiagonal and $B$
tridiagonal. We solve the eigenvalue problem (by means of the inverse power method
with adapted shifts) in order to simply get the eigenvalue with the largest real
part, that is to say, the most unstable one. Precision is improved through
extrapolations on the mesh size $2/(n+1)$. We have obtained the known results
reported by other authors, e.g. \citet{Orszag71}, with an analogous accuracy. The
neutral stability curve, where $\lambda_r = 0$, is presented in
\reffig{estab-lineal}. In this figure, each point in the $Re$--$\alpha$ plane
represents a perturbation of the laminar solution whose stability is decided upon
its position: green points are stable, red ones unstable and blue ones neutrally
stable. We also observe the critical Reynolds number, $Re_{\cri}=5772.22$ for
$\alpha=1.02056$, so that if $Re<Re_{\cri}$ the laminar solution is linearly stable
for any value of $\alpha$. Likewise, for $\alpha\gtrsim1.1$ the laminar flow is
linearly stable for every $Re$.

\figura{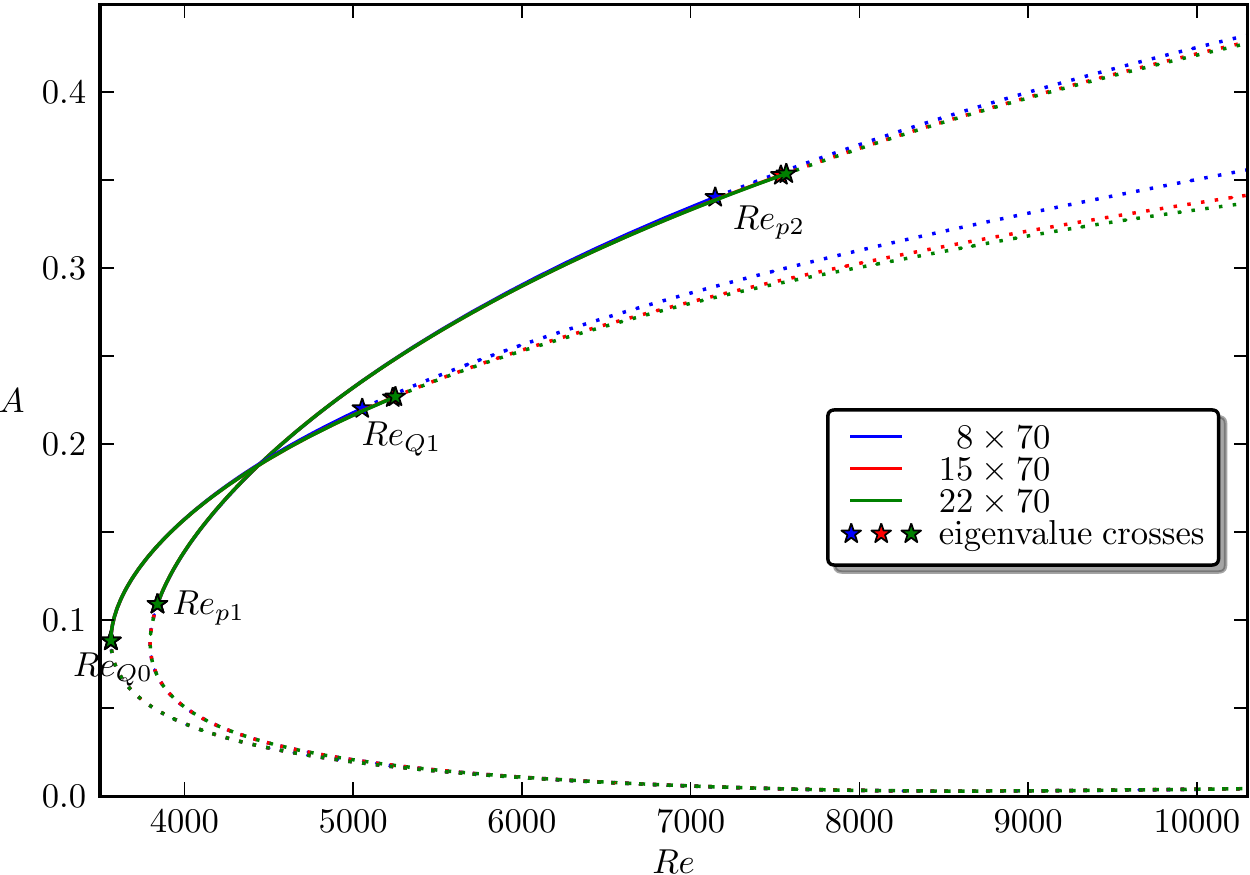}{Bifurcating curve of periodic flows for several
discretizations specified as $N\times M$, $\alpha=1.1$, and based on $Re_p$ and
$Re_Q$. On each curve based on $Re_p$ there are several `*' corresponding to Hopf
bifurcations. They divide the different regions of stability to superharmonic
disturbances, which are also represented in the plot as continuous (stable) and
discontinuous (unstable) lines. In the $Re_Q$ case, the point labeled $Re_{Q0}$
represents a real eigenvalue crossing the imaginary axis, meanwhile $Re_{Q1}$ is a
Hopf bifurcation. Likewise, at $Re_{Q0}$ it is attained the minimum $Re_Q$}

\figura{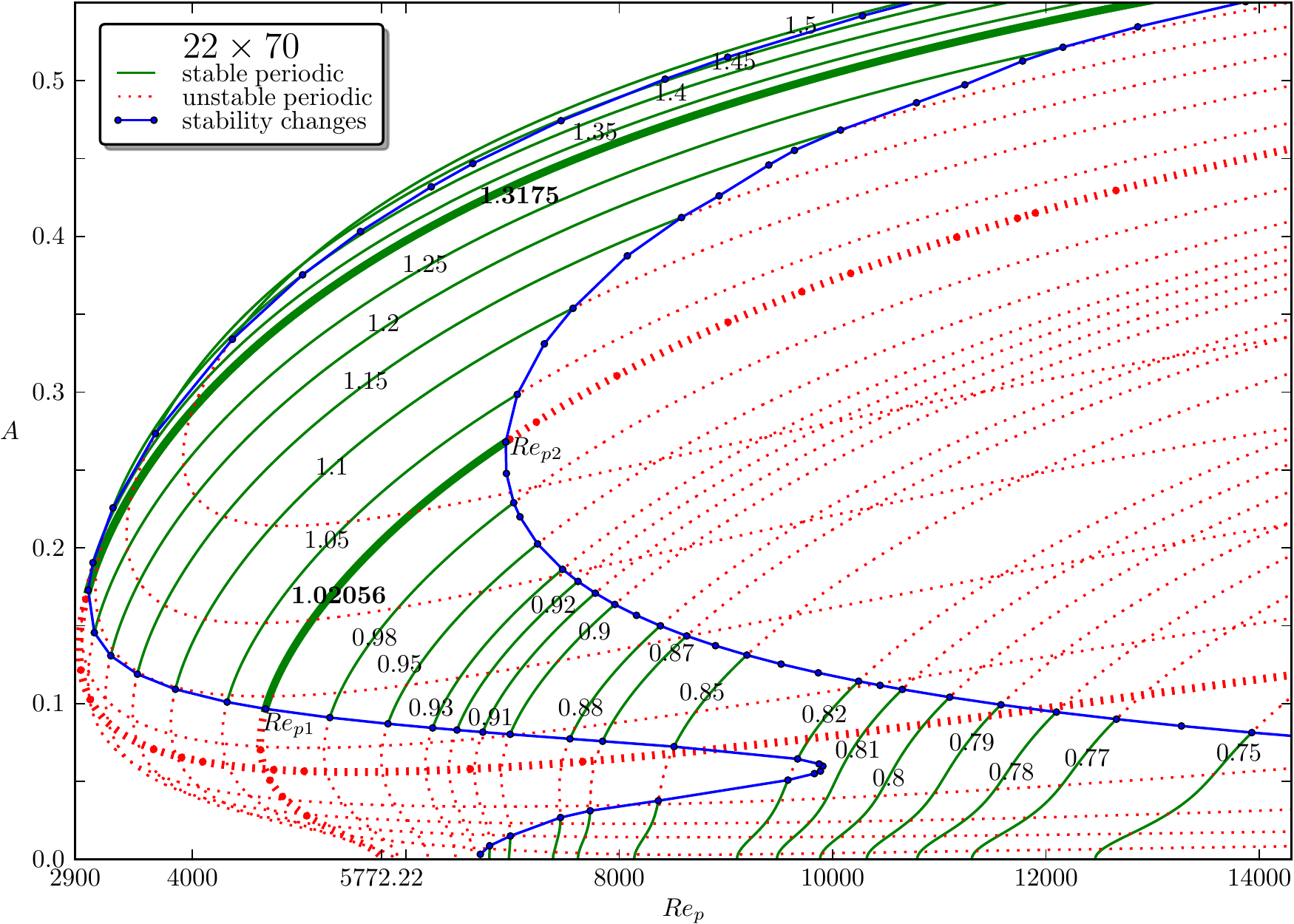}{Bifurcating curves of periodic flows for $Re_p$ and
several values of $\alpha$ specified on each curve. The number of discretization
points is $N=22, M=70$. The `$\bullet$' on each curve represents a change of
stability. Joining those points we obtain $4$ regions: the corresponding solutions
are unstable on the right and bottom-left regions and stable on the top region. In
the intermediate region there are stable and unstable flows, even for a single
$\alpha$, e.g. $\alpha=1.3175$. Curves for $\alpha=1.02056$ and $\alpha=1.3175$ are
traced in thicker lines. For $\alpha=1.02056$ it is attained the critical Reynolds
number at $Re_p = 5772.22$ and for $\alpha=1.3175$ it is approximately reached a
solution at a minimum $Re_p = 2939$}

\subsection{Continuation of travelling waves}

Next, we use the above results to look for periodic solutions in time. Due to the
translational symmetry of the channel in the stream direction (artificial boundaries
in \eqref{CF}), it is showed in \citet{Rand82} that if we have $u(x,y,t)$ such that
\[u(x,y,t+T)=u(x,y,t),\quad\mbox{for all } x\in\R, \quad y\in[-1,1], \quad t\ge0\]
and some $T>0$ (that is, $u(x,y,t)$ is $T$-periodic in time), then it is a rotating
(or travelling) wave, i.e. \begin{equation}u(x,y,t)=u(x-ct,y,0), \quad\mbox{for }
  c=\frac{L}{ T}.\label{OR}\end{equation} Consequently $u(x,y,t)$ is observed as a
stationary solution in a system of reference moving at speed $c$ as it was
introduced in \eqref{NSc}. The converse is also true, namely every stationary
solution of \eqref{NSc} gives rise to a time-periodic solution as is easily
verified. This fact allows us to search for periodic solutions in time as functions
$u(x,y)$ in a Galilean reference at speed $c$, which solve the stationary version of
\eqref{NSc} or, in its discretized form, the stationary version of \eqref{NSred}
\begin{equation}\begin{cases} 0 = U_0 \\ 0 = (I - P_k T_k)^{-1}(\bar U_k - P_k\bar
V_k), & k=1, \dots, N.
\end{cases}\label{NSredE}\end{equation}
Given a fixed $\alpha$, what we have in \eqref{NSredE} is a zeros search problem for
a system of nonlinear equations of dimension $K$ (defined at the end of
\refsecc{numerical_app}). It can be expressed as $H_p(Re,c,U)=0$. Solutions of
(\ref{NSredE}) are locally unique for each $Re$, except translations in the stream
direction. This is due to the fact that any translation of a rotating wave in the
stream direction, gives rise to the same wave at a different time instant. Indeed,
if $u(x,y,0)$ is the starting position of a rotating wave then, by \eqref{OR}, for
every $\theta\in\R$ we have $u(x-\theta,y,0)=u(x,y,\theta/c),$ and thus
$u(x-\theta,y,0)$ lies in the same orbit as $u(x,y,0)$. In order to achieve
uniqueness, we fix one of the coordinates of $U$, by restricting it to a Poincar\'e
section
\begin{equation}\Sigma_1 = \left\{ U = (\bar u_0, \dots, \bar u_N) \ |\ \Re(\bar
  u_{11}) = s_1 \right\}, \label{sec_Poinca_1}\end{equation} ($\Re(u), \Im(u)$ stand
for the real and imaginary part of $u$), for $s_1 \in\R$, a fixed value. We mainly
set $s_1=0$, since this choice gives well conditioned systems. If we fix $Re$ in
\eqref{NSredE}, the number of unknowns, $K$, is the same as the number of equations.
We look for its solutions by means of pseudo-Newton's method, in which we factorize
the resulting linear system using a direct $LU$ decomposition. The first
approximation of the Jacobian matrix $DH_p$ is implemented by finite differences
with extrapolation. Every column of this matrix is obtained evaluating $H_p$ in
parallel in a Beowulf system. Subsequent updates of $DH_p$ are carried out by
Broyden's `good' formula. As a consequence, at each pseudo-Newton step we only have
to apply rank-one updates to the $LU$ factorization of $DH_p$. These improvements
mean an enormous increase in the speed of computations.

If we use $Re$ as a continuation parameter, we can trace the one-parameter curve
$H_p(Re,c,U)=0$. This is implemented numerically by pseudo-arclength continuation.
To this end, we compute the unit-norm tangent vector to the curve. This vector and
previously computed points on the curve, are used to predict the next solution
point, which is finally corrected by pseudo-Newton iterations.

The starting point of those iterations is a numerically integrated periodic
solution. We obtain it using the numerical integrator \eqref{IN}. The initial
condition is taken as a small perturbation of the laminar flow. At the same time we
check our previous results and those reported in \citet{Orszag71} about $Re_{\cri}$.
Namely, taking for example $\alpha=1.02056, Re>Re_{\cri}$, we observe how the
laminar flow is not stable, as it evolves to another steady flow which turns out to
be $T$-periodic in time for some $T$. Up to errors of order $O((\Delta t)^2)$ (those
of the time discretization \eqref{IN}), this $T$-periodic solution satisfies
\eqref{NSredE} and it is thus a rotating wave. We choose it as initial approximation
to a point on the curve $H_p(Re,c,U)=0$.

\begin{table}
\tabcolsep0,25cm\arrayrulewidth0,2pt
\caption{\label{min_Re_per} For some values of $\alpha$, this table shows the
minimum values of $Re_p$ and $Re_Q$ for which exists periodic flow. Calculations
made for $N = 22$ and $M = 70$. The minimum $Re$ attained is marked with `*'.
Results are in good agreement with those reported by \citet{Herbert76} for $Re_p$.}
\begin{tabular}{cc@{\hspace{1cm}}cc}\hline\hline
$\alpha$ &$\min Re_Q$ &$\alpha$ &$\min Re_p$ \\
\hline
1.1000  &3564.5164  &1.1000  &3797.0331 \\
1.2236  &2845.5884  &1.2400  &3048.0073 \\
1.3000  &2647.6068  &1.3092  &2939.3711 \\
1.3424  &2608.9990  &1.3145  &2939.2069 \\
1.3520  &2607.5519  &1.3174  &2939.0345 \\
1.3521  &\omit$\phantom{\kern0,15cml}2607.5516{}^*$\hidewidth  &1.3175
  &\omit$\phantom{\kern0,15cml}2939.0343^*$\hidewidth \\
1.3523  &2607.5520  &1.3177  &2939.0350 \\
1.3534  &2607.5753  &1.3265  &2940.6307 \\
1.5000  &3018.3031  &1.4665  &3526.0725 \\
\hline\hline\end{tabular}\end{table}

\figura{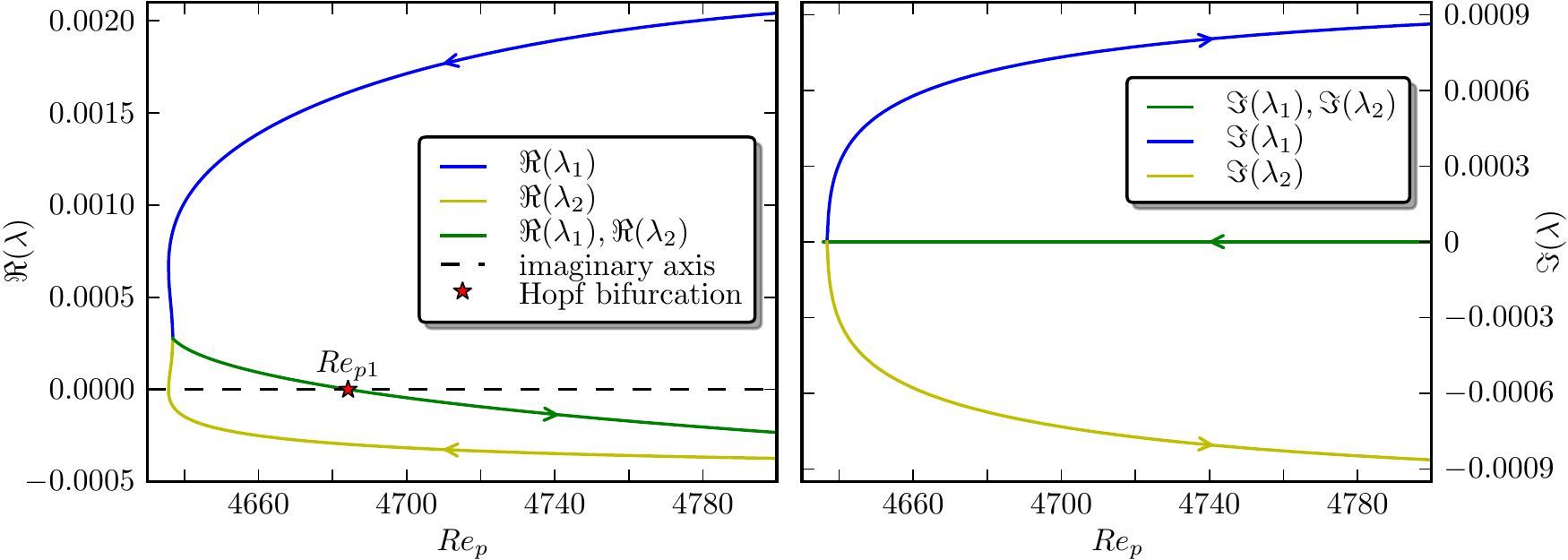}{Real and imaginary part of the two most unstable eigenvalues
($\lambda_1$ and $\lambda_2$) of periodic flows for $N = 22$, $M = 70$, $\alpha =
1.02056$ and $Re_p$ (cf. \reffig{bifur-perp_22x70}). Close to the minimum $Re_p$ of
the amplitude curve, $\lambda_1$ and $\lambda_2$ are unstable and real, and give
rise to a pair of complex conjugate eigenvalues. The first crossing through the
imaginary axis is associated to a single real eigenvalue ($\lambda_2$), meanwhile
the second one corresponds to the $Re_{p1}$ Hopf bifurcation ($\lambda_1$ and
$\lambda_2$). Arrows for both graphs, point to the direction of increasing amplitude
in \reffig{bifur-perp_22x70}}

Given a profile of velocities $(u,v)$ we define its amplitude $A$, as the distance
to the laminar profile $(u_b,0)$ in the $L^2$-norm
\begin{equation}A=\frac{1}{2L}\|(u-u_b,v)\|_2, \qquad \|(u,v)\|_2^2=
\int_0^L\int_{-1}^1\left[u(x,y)^2 +
v(x,y)^2\right]\MR{d}y\,\MR{d}x.\label{amplitud}\end{equation}
For a rotating wave as defined in \eqref{OR}, its amplitude does not depend on
time since for fixed $t$:
\begin{equation}\begin{array}{r@{}l}
\DS\int_0^L\int_{-1}^1 u(x,y,t)^2 \MR dy\,\MR dx
\buildrel\mbox{\scriptsize1}\over=&\DS \int_0^L\int_{-1}^1 u(x-ct,y,0)^2\MR dy\, \MR
dx \nonumber \\ [3ex] \DS\buildrel\mbox{\scriptsize2}\over=
&\DS\int_{-ct}^{L-ct}\int_{-1}^1 u(\tilde x,y,0)^2 \MR{d}y\,\MR{d}\tilde
x\buildrel\mbox{\scriptsize3}\over= \int_0^L\int_{-1}^1 u(\tilde x,y,0)^2
\MR{d}y\,\MR{d}\tilde x.\end{array} \end{equation} In step 1 we
apply definition \eqref{OR}. For step 2 we make the change of variable $\tilde
x=x-ct$, and because $u$ is $L$-periodic in $x$ we have step 3.

\figura{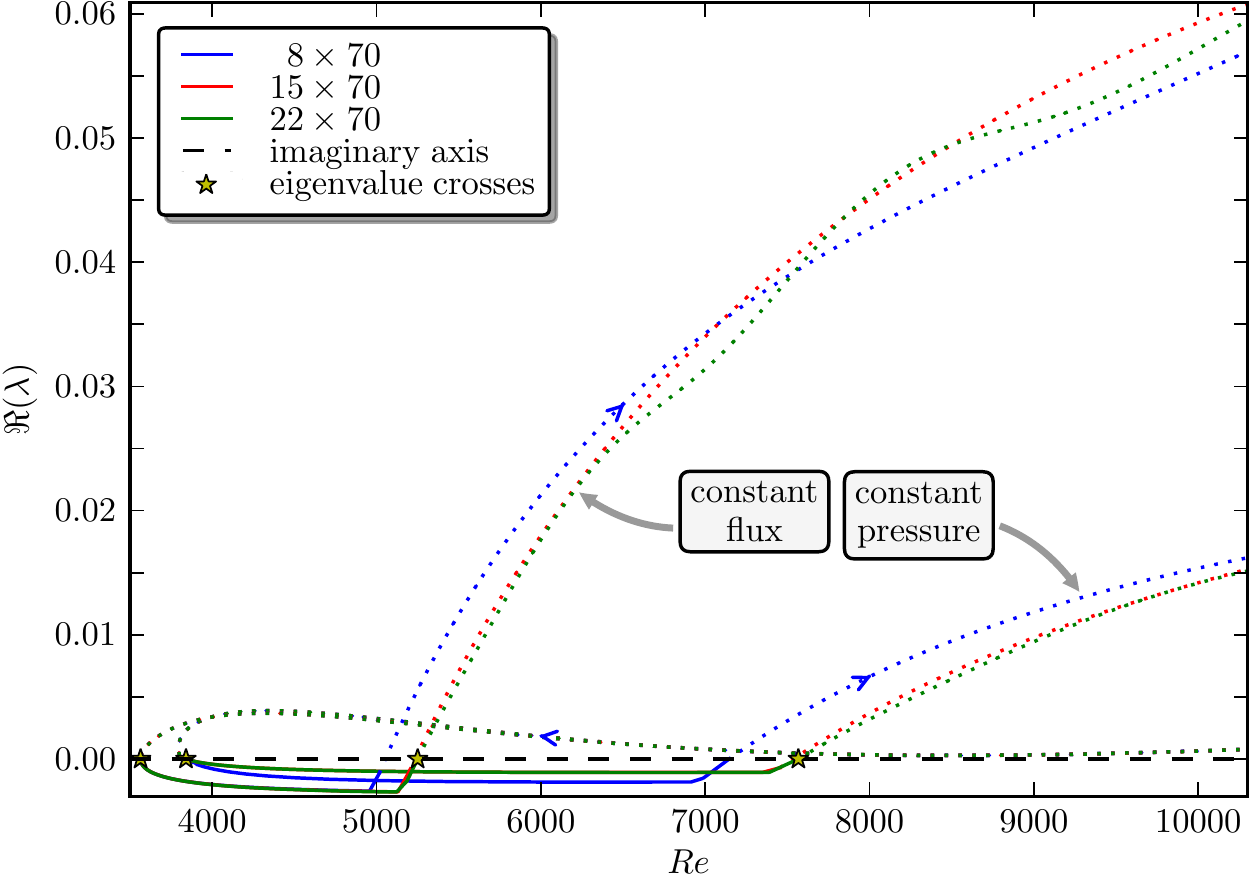}{Real part of the most unstable eigenvalue for periodic flows
for $\alpha = 1.1$ and $N \times M$ as specified. For the $Re_p$ case, the two
crossings of each curve with the imaginary axis corresponds to the first two `*' of
\reffig{bifur-per_1_1}, and are detailed in \reftab{bifur_Hopf} as $Re_{p1}$ and
$Re_{p2}$. In the analogous plot based on $Re_Q$, the two crossings of each curve
with the imaginary axis corresponds to the first two `*' of \reffig{bifur-per_1_1}
and are also specified in \reftab{bifur_Hopf} as $Re_{Q0}$ and $Re_{Q1}$. Continuous
and discontinuous lines refer respectively to stable and unstable periodic solutions
associated to each point $(Re,\Re(\lambda))$. Arrows point to the direction of
increasing amplitude in \reffig{bifur-per_1_1}}

\subsection{Stability of periodic solutions}\label{estab_sol_per}

We notice that zeros of system \eqref{NSredE} can either correspond to stable or
unstable time-periodic solutions. With a stable solution it is meant the one for
which any small disturbance ultimately decays to zero, whereas if some of those
disturbances remain permanently away from zero, it is called unstable.

\figura{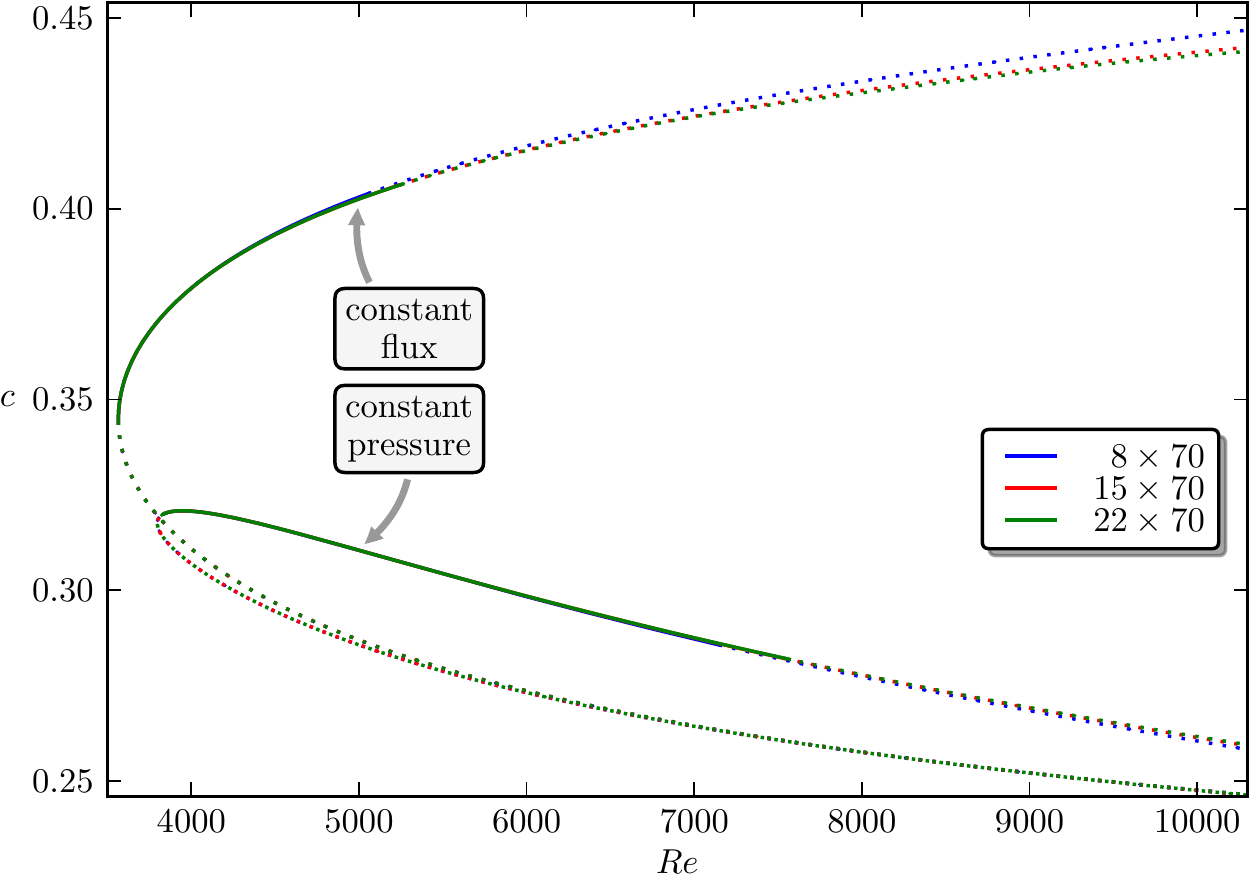}{Speed of the observer $c$, for parameters and periodic solutions
as in \reffig{bifur-per_1_1}. On both $Re_p$ and $Re_Q$ cases, upper and lower
branches correspond to the respective ones of amplitudes. Unlike $Re_p$, for the
constant flux case the upper branch increases with $Re_Q$. Continuous and
discontinuous lines refer respectively to stable and unstable periodic solutions
associated to each $c$}

To decide whether a time-periodic flow $u$ is stable or not, we consider it as a
steady solution for its appropriate $c=L/T$ and obtain the eigenvalues of its
Jacobian matrix. This matrix is computed analytically linearizing \eqref{NSred}
around $u$. If every eigenvalue has negative real part, the periodic flow is stable
to disturbances of the same wavenumber $\alpha$ but, if there is an eigenvalue with
positive real part, the solution is unstable. Let us mention that there is always a
zero eigenvalue which corresponds to the lack of uniqueness of the time-periodic
flow due to translations. Setting $\alpha = 1.1$, the bifurcating diagram for the
periodic flows in the $Re$-$A$ plane together with the stability changes are
represented in figure \ref{bifur-per_1_1} for both formulations in terms of $Re_p$
and $Re_Q$. Due to relations \eqref{RelRey} and \eqref{RelVel} we only need to
compute travelling waves $u^p(x,y,t)$ at speed $c$ for $Re_p$, since \eqref{RelVel}
gives the associated $u^Q(x,y,t)$, which is easily checked to be a travelling wave
at speed $rc$ for $Re_Q$, being $r = Re_p/Re_Q$. As well as computing the
eigenvalues, we confirm the stability of a periodic flow using the numerical
integrator \eqref{IN}. A more complete picture of the different connections among
stable and unstable solutions is given in \citet{Casas04}.

By simply taking a known travelling wave for some $\alpha$ as initial guess and
moving slightly $\alpha$, we can find periodic solutions for different values of
$\alpha$. These are shown in \reffig{bifur-perp_22x70}, together with their
stability. In addition, in \reftab{min_Re_per} we have computed, for several values
of $\alpha$, the corresponding minimum value of $Re_p$ and $Re_Q$ (denoted as
$Re_{min}(\alpha)$) along the amplitude curves (see \reffig{bifur-perp_22x70}). In
turn, $Re_{min}(\alpha)$, is minimized as a function of $\alpha$. In this way we
obtain the absolute minimum $Re_p$ and $Re_Q$ for which there exists periodic
solution. These minimum values are marked with `*' in \reftab{min_Re_per}. For
$Re_p$, \citet{Herbert76} obtained the minimum value at $Re_p = 2934.80$ for $\alpha
= 1.3231$ and $N \times M = 4 \times 40$ as the spectral spatial discretization for
the stream function: this discretization is analogous to the one used in the present
work. We observe that our value of $Re_p$ differ from Herbert's not more than $0.15
\%$.

Previously, \citet{Soibelman91} found similar bifurcations of travelling waves for
$\alpha = 1.1$, and the critical Reynolds number for which there are time-periodic
solutions: $Re_p \approx 2900$ for $\alpha \approx 1.3$, and $Re_Q \approx 2600$. We
remark that in \reffig{bifur-perp_22x70} for $\alpha=1.3175$ there exists an
attracting periodic solution for $Re_p = 3024$ which corresponds to $Re_Q = 2630$.
For the case of the laminar flow, we encounter the classical results of
\citet{Orszag71} about the critical Reynolds being at $Re_{\cri} = 5772.22$ for
$\alpha=1.02056$. On the one hand, we observe in \reffig{bifur-perp_22x70} that the
bifurcation curve of periodic flows reaches the laminar solution at the above
mentioned $Re_{\cri}$ and in addition, the laminar solution is checked to be stable
when $Re<Re_{\cri}$ and unstable if $Re>Re_{\cri}$. This Hopf bifurcation at
$Re_{\cri}$ is called subcritical, because the branch of periodic solutions
emanating at it decreases in $Re$. When the new branch increases in $Re$, we call it
supercritical bifurcation. We also notice that for $\alpha=1.02056$, it is reached
the minimum Reynolds number where the transition from stable to unstable laminar
flow takes place, what was formerly presented in \reffig{estab-lineal}. For
$\alpha\gtrsim 1.1$ the curve of periodic solutions does not reach the laminar flow.
This is in agreement with the situation shown in \reffig{estab-lineal}, since for
$\alpha \gtrsim 1.1$ the laminar flow is linearly stable for every $Re$. In
\reffig{bifur-perp_22x70} we check as well that for $\alpha \gtrsim 0.91$ the curve
of periodic flows bifurcates subcritically from the laminar flow, but for $\alpha
\lesssim 0.91$ the Hopf bifurcation is changed into supercritical. Precisely for
$\alpha\approx 0.915$ it is born a new change of stability on the curve of
travelling waves at $Re_p \approx 6700$. The behaviour around this new point is
saddle-node bifurcation, analogous to $Re_{Q0}$, and will be described in the
following subsection.

\begin{table}
\arrayrulewidth0,2pt\tabcolsep0,15cm
\caption{\label{bifur_Hopf}Minimum Reynolds number $Re_{Q0}$ and the three Hopf
bifurcations of periodic flows at $Re_{p1}$, $Re_{Q1}$, $Re_{p2}$, together with
associated parameters $c$ and $\tau$ for $\alpha = 1.02056, 1.1$, $M = 70$ and
several $N$. The values reported in \citet{Soibelman91} for $M=70$ are also
included.}
\begin{tabular}{r@{\hspace{7ex}}cc@{\hspace{7ex}}ccc@{\hspace{7ex}}ccc@{\hspace{7ex}}ccc}\hline\hline
$N$ &$Re_{Q0}$ &$c_0$ &$Re_{p1}$ &$c_1$ &$\tau_1$ &$Re_{Q1}$ &$c_1$ &$\tau_1$ &$Re_{p2}$
&$c_2$ &$\tau_2$ \\ [1ex]

\multicolumn{12}{c}{Present work for $\alpha=1.02056$} \\ \hline
 4 & 4439.1 & 0.31544 & 4701.7 & 0.29729 & 9662.43 & 5620.0 & 0.37239 & 13.52 & 7450.1 & 0.28091 & 17.93 \\
 5 & 4387.1 & 0.31934 & 4684.6 & 0.29881 & 8983.08 & 5108.0 & 0.36054 & 14.00 & 6347.5 & 0.29014 & 17.39 \\
 6 & 4393.9 & 0.31890 & 4686.1 & 0.29868 & 9080.07 & 5085.8 & 0.36257 & 13.78 & 6355.4 & 0.29013 & 17.22 \\
 7 & 4396.3 & 0.31845 & 4684.1 & 0.29853 & 9143.27 & 5243.2 & 0.36468 & 13.68 & 6643.2 & 0.28783 & 17.33 \\
 8 & 4395.4 & 0.31854 & 4684.1 & 0.29857 & 9125.21 & 5335.6 & 0.36610 & 13.67 & 6827.8 & 0.28608 & 17.49 \\
 9 & 4395.2 & 0.31858 & 4684.3 & 0.29858 & 9120.06 & 5341.9 & 0.36637 & 13.70 & 6845.0 & 0.28592 & 17.55 \\
10 & 4395.2 & 0.31857 & 4684.2 & 0.29858 & 9120.76 & 5353.9 & 0.36653 & 13.68 & 6865.9 & 0.28581 & 17.54 \\
11 & 4395.2 & 0.31858 & 4684.2 & 0.29858 & 9120.24 & 5371.9 & 0.36680 & 13.64 & 6900.7 & 0.28555 & 17.52 \\
12 & 4395.2 & 0.31858 & 4684.2 & 0.29858 & 9119.98 & 5385.0 & 0.36703 & 13.62 & 6926.8 & 0.28534 & 17.52 \\
\multicolumn{1}{@{\hspace{2ex}}l}{\vdots}&&&&&     &\vdots  &\vdots   &\vdots &\vdots  &\vdots   &\vdots \\
18 &        &         &        &         &         & 5253.3 & 0.40634 & 11.78 & 6936.0 & 0.28530 & 17.51 \\
19 &        &         &        &         &         & 5389.8 & 0.36711 & 13.61 & 6934.9 & 0.28532 & 17.51 \\
20 &        &         &        &         &         & 5390.4 & 0.36713 & 13.61 & 6936.3 & 0.28530 & 17.51 \\
21 &        &         &        &         &         & 5390.3 & 0.36712 & 13.61 & 6935.9 & 0.28531 & 17.51 \\
22 &        &         &        &         &         & 5390.3 & 0.36712 & 13.61 & 6936.1 & 0.28531 & 17.51 \\
[1ex]

\multicolumn{12}{c}{Present work for $\alpha=1.1$} \\ \hline
 4 & 3603.5 & 0.34248 & 3864.5 & 0.31927 & 7429.40 & 5812.6 & 0.41532 & 11.62 & 8946.1 & 0.26985 & 17.88 \\
 5 & 3564.7 & 0.34467 & 3841.0 & 0.32004 & 7090.00 & 5296.6 & 0.40633 & 11.65 & 7670.3 & 0.28059 & 16.87 \\
 6 & 3562.8 & 0.34506 & 3841.4 & 0.32018 & 7065.52 & 4840.4 & 0.40235 & 11.93 & 6732.7 & 0.28927 & 16.59 \\
 7 & 3564.9 & 0.34475 & 3840.8 & 0.32010 & 7100.29 & 4905.9 & 0.40270 & 11.88 & 6844.8 & 0.28863 & 16.58 \\
 8 & 3564.7 & 0.34474 & 3840.6 & 0.32010 & 7098.78 & 5054.0 & 0.40392 & 11.85 & 7145.8 & 0.28568 & 16.75 \\
 9 & 3564.5 & 0.34477 & 3840.6 & 0.32010 & 7096.05 & 5107.8 & 0.40457 & 11.86 & 7264.9 & 0.28445 & 16.87 \\
10 & 3564.5 & 0.34477 & 3840.6 & 0.32011 & 7095.79 & 5120.9 & 0.40477 & 11.87 & 7292.5 & 0.28424 & 16.90 \\
11 & 3564.5 & 0.34477 & 3840.6 & 0.32011 & 7095.68 & 5157.2 & 0.40519 & 11.84 & 7368.9 & 0.28358 & 16.92 \\
12 & 3564.5 & 0.34477 & 3840.6 & 0.32011 & 7095.56 & 5189.1 & 0.40555 & 11.82 & 7436.6 & 0.28298 & 16.94 \\
\multicolumn{1}{@{\hspace{2ex}}l}{\vdots}&&&&&     &\vdots  &\vdots   &\vdots &\vdots  &\vdots   &\vdots \\
21 &        &         &        &         &         & 5246.8 & 0.40624 & 11.78 & 7558.9 & 0.28197 & 16.97 \\
22 &        &         &        &         &         & 5250.4 & 0.40629 & 11.78 & 7567.0 & 0.28190 & 16.98 \\
23 &        &         &        &         &         & 5249.5 & 0.40628 & 11.78 & 7565.0 & 0.28192 & 16.97 \\
24 &        &         &        &         &         & 5250.6 & 0.40630 & 11.78 & 7567.4 & 0.28190 & 16.98 \\
25 &        &         &        &         &         & 5249.6 & 0.40628 & 11.78 & 7565.2 & 0.28192 & 16.98 \\
[1ex]

\multicolumn{12}{c}{\citet{Soibelman91} for $\alpha=1.1$} \\ \hline
2 &&&\multicolumn{2}{@{}l}{3630} & 4742.32 & \multicolumn{2}{@{}l}{5600} &
\multicolumn{1}{@{\hspace{1ex}}l}{20.6} &\multicolumn{2}{@{}l}{9400} &35.50 \\
3 &&&\multicolumn{2}{@{}l}{3800} & 4935.43 & \multicolumn{2}{@{}l}{6250} &
\multicolumn{1}{@{\hspace{1ex}}l}{12.5} &\multicolumn{2}{@{}l}{9675} &17.65 \\
4 &&&\multicolumn{2}{@{}l}{3775} & 4875.63 & \multicolumn{2}{@{}l}{5875} &
\multicolumn{1}{@{\hspace{1ex}}l}{13.4} &\multicolumn{2}{@{}l}{9592} &16.54 \\ [1ex]
\hline\hline\end{tabular}\end{table}

\subsection{Hopf bifurcations}

Now let us discuss the bifurcation diagram shown in \reffig{bifur-perp_22x70} for
$\alpha=1.02056$. First, the laminar solution becomes unstable at the critical value
$Re_{\cri} = 5772.22$, due to a Hopf bifurcation that gives rise to a unstable
family of periodic orbits. This family continues backwards (with respect to $Re$,
i.e. it is subcritical) until $Re_p\approx 4636$, where a turning point is reached
and $\Re(\lambda_2)$ (described in \reffig{vap-en-Hopf}) crosses the imaginary axis.
Before arriving at this turning point, there is a single eigenvalue, $\lambda_1$, on
the real positive axis, while the remaining ones have negative real part (we ignore
the eigenvalue at $0$ arising from the lack of uniqueness of periodic flows). On
traversing through the turning point, a real and negative eigenvalue ($\lambda_2$)
becomes real positive, so the number of unstable eigenvalues is now two. Shortly
after that, these two unstable eigenvalues collide and become a conjugate complex
pair (still with positive real part), and then they cross the imaginary axis for
$Re_p \approx 4684$ producing a new Hopf bifurcation at the point $Re_{p1}$ on
\reffig{bifur-perp_22x70}. Between $Re_{p1}$ and $Re_{p2} \approx 6936$, the family
of periodic orbits is stable to disturbances of the same wavelength. At $Re_{p2}$,
there is another Hopf bifurcation produced by a pair of conjugate eigenvalues
crossing the imaginary axis. These bifurcations persist, as shown in
\reftab{bifur_Hopf}, when $M,N$ are increased and no new ones seem to appear in this
range.

The case of constant flux is qualitatively different. For $Re_Q$ the bifurcating
diagram of periodic solutions has a turning point at a minimum value of $Re_Q$,
which we designate as $Re_{Q0}$ (cf. \reffig{bifur-per_1_1}). The lower branch of
periodic solutions is unstable with only one unstable real eigenvalue and the upper
is initially stable, being also real the most unstable eigenvalue. On traversing the
bifurcating curve towards the upper branch, this real positive eigenvalue becomes
negative at the turning point. The upper branch is kept stable until a subsequent
Hopf bifurcation appears at certain value $Re_{Q1}$. For $Re_Q \gtrsim 7000$ and
$\alpha=1.1$, we have detected more Hopf bifurcations which we do not consider in
this study. However for the range included in \reffig{bifur-per_1_1} all periodic
flows for $Re_Q>Re_{Q1}$ are unstable.

\citet{Pugh88} pointed out that the null eigenvalue at $Re_{Q0}$ has algebraic
multiplicity $2$ and geometric multiplicity $1$. We can consider that eigenvalue
simple (with algebraic and geometric multiplicity 1) if we ignore the constant zero
eigenvalue due to a trivial phase shift of the flow in the stream direction. The
suppression of this trivial null eigenvalue can be made by restricting equations
\eqref{NSred} to the closed linear manifold $\Sigma_1$ defined as a Poincar\'e
section in \eqref{sec_Poinca_1}. According to bifurcation theory, at a simple
eigenvalue like this one we have no equilibrium point for $Re_Q < Re_{Q0}$ and two
equilibrium points for $Re_Q > Re_{Q0}$: this situation corresponds to a saddle-node
bifurcation and no new branches of solutions come out from $Re_{Q0}$.

In \reftab{bifur_Hopf} are shown $Re_p$, $Re_Q$, the speed of the observer $c$ and
the period $\tau$ of the bifurcated solution corresponding to the first Hopf
bifurcations for several values of $N$ and $M = 70$. Taking $M = 70$, Chebyshev
modes seem to be enough to attain convergence in the results. The values obtained by
\citet{Soibelman91} are also presented for comparison. For $\alpha=1.1$ we observe
convergence of our results on the different Hopf bifurcations considered as $N$ is
increased. In all cases there are substantial differences with
\citeauthor{Soibelman91}'s \cite{Soibelman91} results, being in more agreement for
the lowest $Re$. We remark the slow convergence of the Fourier series to the
bifurcation values as $N$ is increased. At the same time we have also obtained
convergence in the qualitative behaviour: the subcritical or supercritical character
of all the studied Hopf bifurcations remain unaltered as $M,N$ are increased.

Formulas \eqref{RelRey} and \eqref{RelVel} provide again the correspondence between
bifurcation points at $Re_Q$ and $Re_p$ (cf. \reftab{bifur_Hopf}). For instance at
$Re_{p1}=3840.6$ for $N = 12$, $M = 70$ and $\alpha = 1.1$, the periodic solution
transformed by these formulas furnish a periodic solution at $Re_Q = 3564.5$ and $c
= 0.34490$, values in good agreement with $Re_{Q0}$ reported in \reftab{bifur_Hopf}.
Likewise, the transformed periodic solution for $N=25$, $\alpha=1.1$ at
$Re_{Q1}=5249.6$ gives rise to $Re_p = 7565.5$ and $c = 0.28192$, again in good
precision with respect to $Re_{p2}$.

The different stability changes marked as a blue dot in \reffig{bifur-perp_22x70},
roughly divide the $Re_p$-$A$ plane in four regions. On the right and bottom-left
regions the corresponding periodic solutions are unstable (filled with red points),
meanwhile they are stable on the top region (only green points). In the intermediate
region there are both stable and unstable flows, even for a single $\alpha$, e.g.
$\alpha=1.3175$. Through this classification, given a periodic flow with its
associated $(Re_p,A)$, we can deduce its stability, independently on $\alpha$ in
some cases. On traversing the $Re_{p1}$ blue curve in the direction of increasing
amplitudes, up to the relative maximum on that curve attained at $Re_{p1} \approx
9909$ for $\alpha \approx 0.814$, we encounter saddle-node bifurcations for
$\alpha\in [0.88,0.915]$ at a relative maximum on each bifurcation curve. For
$\alpha\in [0.714,0.88]$, the former relative maximum disappears and the saddle-node
bifurcation turns into a Hopf one. The rest of the $Re_{p1}$ curve is made up of
Hopf bifurcations for the studied values $\alpha\le 1.7$. The minimum $Re_{p1}
\approx 3024$ is reached precisely for $\alpha \approx 1.3175$, where the minimum
periodic flow was found in \refsecc{estab_sol_per}. The $Re_{p2}$ blue curve is only
constituted of Hopf bifurcations for the $\alpha\in[0.74,1.3175]$ considered. The
minimum $Re_{p2} \approx 6936$ is achieved again for the critical $\alpha \approx
1.02056$.

The maximum growth rate (real part of the most unstable eigenvalue) for each
periodic flow is presented in \reffig{vap-per_1_1} for the same parameters as in
\reffig{bifur-per_1_1}. For the most unstable eigenvalue $\lambda$, $\Re(\lambda)$
crosses the imaginary axis twice, on the values $Re_{p1}$ and $Re_{p2}$ for $Re_p$
and on $Re_{Q0}$ and $Re_{Q1}$ for $Re_Q$. Those diagrams represent the degree of
instability of each flow. Comparing to \reffig{bifur-per_1_1}, we observe that at
the same $Re$ on the upper branch of amplitudes, periodic solutions based on $Re_Q$
are more unstable than the associated ones based on $Re_p$. On the other hand, on
the lower branch of amplitudes, at the same $Re$, both curves of $\Re(\lambda)$
visually coincides for $Re \gtrsim 5000$. This behaviour is also reflected in
figures \ref{bifur-per_1_1} and \ref{c-per_1_1}. In this last figure we present
qualitatively different curves for the speed $c$ in $Re_p$ and $Re_Q$ cases. In the
first case. $c$ is decreasing in both branches of solutions in
\reffig{bifur-per_1_1}. However, for $Re_Q$ the shape of the $c$-curve is similar as
the $A$-curve in \reffig{bifur-per_1_1}.

\section{Quasi-periodic solutions}\label{solcper}

In this section, we study the quasi-periodic flows that appear at the Hopf
bifurcations of rotating waves shown in \refsecc{solper}. They are found as
time-periodic orbits in an appropriate Galilean reference, which simplifies
enormously their search. Those time-periodic orbits are obtained as fixed points of
a Poincar\'e section, by means of a pseudo-Newton method. We have traversed
bifurcating branches of quasi-periodic solutions for the Hopf bifurcations at
$Re_{p1}$, $Re_{p2}$ and $Re_{Q1}$ defined in \refsecc{solper}. We have obtained
different qualitative results than the ones reported in \citet{Soibelman91}. For
$\alpha=1.1$ and $Re_p$ they found that the bifurcation at $Re_{p1}$ to
quasi-periodic solutions is subcritical: they obtained quasi-periodic solutions for
$Re_p$ before the bifurcation point. In consequence, close to that point, those
bifurcated flows are stable. In the present study, by increasing the number of
Fourier modes $N$, we have achieved a supercritical Hopf bifurcation at $Re_{p1}$:
the bifurcating quasi-periodic flows are located for $Re_p$ after the bifurcation
point and therefore close to it they are unstable. This is treated in
\refsecc{bif_Re_p1}. For the second Hopf bifurcation at $Re_{p2}$, in agreement with
\citet{Soibelman91}, the quasi-periodic orbits are found for $Re_p$ greater than the
bifurcation point. More details are given in \refsecc{bif_Re_p2}. The behaviour at
$Re_{Q1}$ (considered in \refsecc{bif_Re_Q1}) is analogous to that of $Re_{p2}$.
However for $Re_Q > Re_{Q1}$ large enough we have detected another Hopf bifurcation
to tori with $3$ basic frequencies. The stability of quasi-periodic flows to
superharmonic disturbances is estimated by means of the linear part of the
Poincar\'e map and also with a full numerical simulation of the fluid.

\subsection{Reduction to periodic and numerical procedures}

We use again the spatio-temporal symmetry of our system, due to the artificial
boundaries of the channel. Considering this symmetry in \citet{Rand82} it is proved
that every solution $u(x,y,t)$ that lies on an isolated invariant $2$-torus (a
quasi-periodic solution), not asymptotic to a rotating wave, is a modulated wave,
that is to say, there exists $\tau>0$ and $\phi \in \R$ such that
\begin{equation}u(x,y,n\tau+t)=u(x-n\phi,y,t) \quad\mbox{for every }
n\in\Z.\label{OM}\end{equation} Hence, this kind of wave has the property that, may
be viewed as a $\tau$-periodic wave in time, in a frame of reference moving at speed
$c = (pL+\phi) / \tau$, for any integer $p$. In effect, defining $\tilde x = x-ct$
for that value of $c$, and $\tilde u(\tilde x,y,t) = u(\tilde x+ct,y,t)$ as the
velocity in the moving frame of reference at speed $c$, it turns out that
\begin{equation}\tilde u(\tilde x,y,\tau)\buildrel\mbox{\scriptsize1}\over=
u(\tilde x+c\tau,y,\tau)\buildrel\mbox{\scriptsize2}\over= u(\tilde x + pL +
\phi,y,\tau)\buildrel\mbox{\scriptsize3}\over= u(\tilde
x,y,0)\buildrel\mbox{\scriptsize4}\over= \tilde u(\tilde x,y,0).
\label{RedaPer}\end{equation} In steps 1, 4 we use the previous definition of
$\tilde u$. Substituting the previously defined $c$ we obtain step 2, and because
$u$ is $L$-periodic in $x$ and a modulated wave, one gets step 3. Consequently we
have proved that $\tilde u(\tilde x,y,t)$ is a $\tau$-periodic function of $t$.

\figura{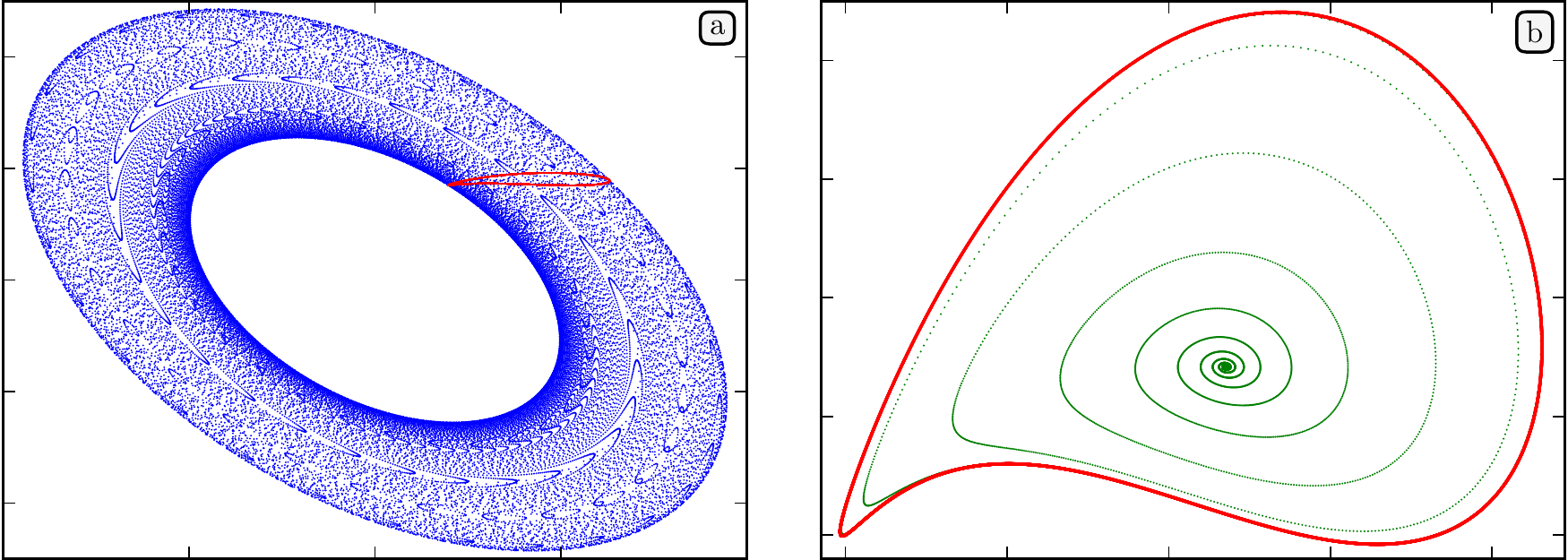}{Two representations of the solution vector $U(t)$ of a
quasi-periodic flow projected on the plane of two selected coordinates namely,
$M-1+(2N-3)(M-2)+3$, and $M-1+2(M-2)+5$, for $Re_p = 3865, \alpha=1.1, N=8, M=70$.
The range of values for (a) is $[-0.002,0.002] \times [-0.001,0.001]$ and
$[0.00037,0.00129] \times [0.000338,0.000385]$ for (b). In (a) a dot is plotted each
$\Delta t=0.02$ time units, meanwhile in (b) only for $t$ such that $U(t) \in
\Sigma_1$. The red closed curve in (b) is obtained integrating repeatedly for $t
\in[0,\tau]$. This curve is also represented as a reference on the right centre of
(a)}

In order to look for periodic flows satisfying \eqref{RedaPer} in a Galilean
reference at speed $c$, we make use of the Poincar\'e section $\Sigma_1$ defined in
(\ref{sec_Poinca_1}). In this case, we only consider points on $\Sigma_1$ when they
cross the section in a particular direction as time increases, namely, from $s_1<0$
to $s_1>0$. Likewise we define the associated Poincar\'e map $P_c : \Sigma_1
\longrightarrow \Sigma_1$ as follows: starting from an initial condition $U=U(0) \in
\Sigma_1$ we integrate \eqref{NSred} for fixed parameters $Re$, $\alpha$ and $c$,
until a time $t_c$ such that $\smash{\tilde U(t_c)} \in \Sigma_1$ ($\smash{\tilde
  U(t_c)}$ represents the evolution of $U(0)$ in a Galilean reference at speed $c$)
for the $n_c$-th time (i.e. after $n_c$ crosses with $\Sigma_1$), where $n_c$ is a
positive integer which represents the minimum number of times needed for the flow to
return close to the initial point $U(0)$ (the meaning of `close' will be specified
in \refsecc{bif_Re_p1}). We then set $P_c(U(0)) = \smash{\tilde U(t_c)}$. In this
way we have reduced the search of quasi-periodic flows to a zeros finding problem
for the map $H_q$ defined as
\begin{equation}0 = H_q(Re,c,U)\buildrel\mbox{\scriptsize def}\over=
P_c(U) - U = \tilde U(t_c) - U(0). \label{ceros}\end{equation}

\figura{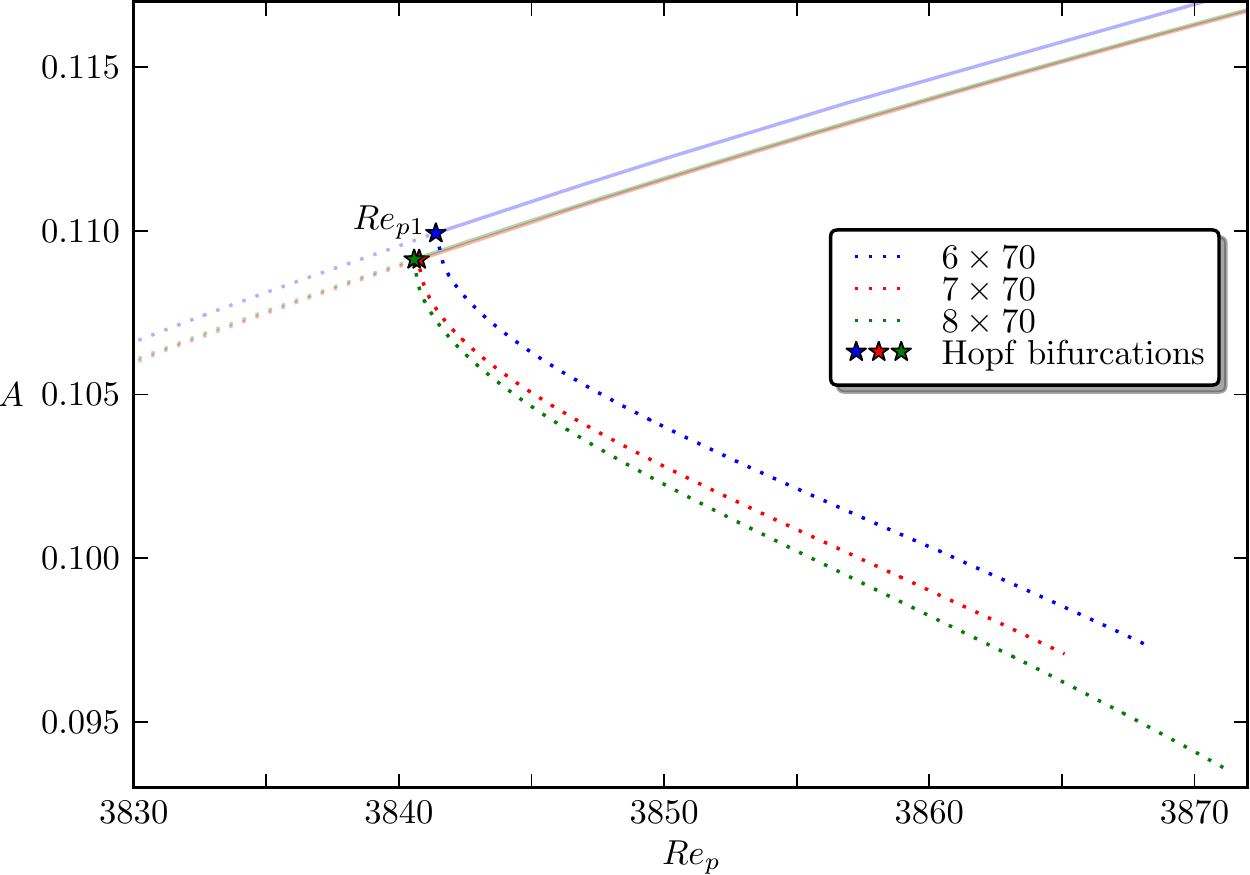}{Bifurcated branches of quasi-periodic flows at the
supercritical Hopf bifurcation of periodic flows at $Re_{p1}$. Each curve represents
$A$ as a function of $Re_p$ and has specified $N \times M$. Both branches of
periodic (in lighter colors) and quasi-periodic flows are presented. Calculations
are shown for $\alpha=1.1$ and $\Delta t=0.02$. The `*' correspond to the Hopf
bifurcation at $Re_{p1}$ presented in \reffig{bifur-per_1_1}. For the range shown,
the bifurcating branch consists of unstable quasi-periodic orbits (dotted lines),
whose amplitude decreases with $Re_p$}

From \eqref{OM} we have that, if $P_c(U(0)) = \tilde U(\tau) = U(0)$, for some
$\tau$, then as well $P_c(V(0)) = \tilde V(\tau) = V(0)$ for $V(0) = P_0^k(U(0))$
and $k$ any positive integer. Indeed, we can express \eqref{RedaPer} as $\tilde
u(\tilde x,y,\tau+t_1) = \tilde u(\tilde x,y,t_1)$ for every $t_1$. From here, since
$V(0) = U(t_1)$ for some $t_1$, we immediately obtain $P_c(V(0)) = V(\tau) = V(0)$.
Therefore we can generate different points on the same orbit as a solution of
\eqref{ceros}. We avoid this lack of uniqueness by restricting $P_c: \Sigma_1 \cap
\Sigma_2 \longrightarrow \Sigma_1$, for $\Sigma_2$ a Poincar\'e section (analogous
to $\Sigma_1$) defined by
\begin{equation}\Sigma_2 = \left\{ U = (\bar u_0, \dots, \bar u_N) \ |\ S = 0
\right\}, \label{sec_Poinca_2}\end{equation} where we set $S = \Re(\bar u_{N,M/2-1})
- s_2$, for $s_2 \in\R$ a suitable quantity. For $Re$ close to the studied Hopf
bifurcations, we have chosen $s_2 = \Re(\bar u^p_{N,M/2-1})$, with $u^p \in
\Sigma_1$ the travelling wave at the exact $Re$ where the bifurcation takes place.
The reason for this choice is merely to preserve continuity in the amplitude
diagrams described next.

\figura{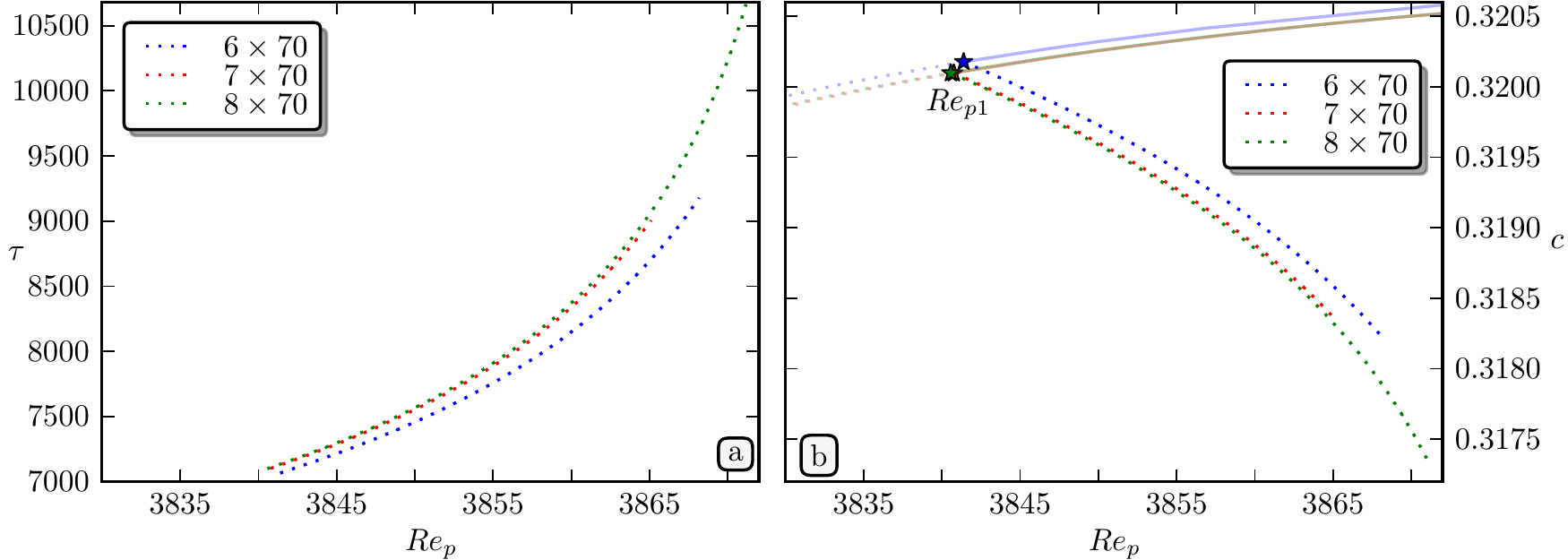}{Different curves related to \reffig{bifur-cper1_1_1} for
$N\times M$ points as specified. (a) $\tau$ represents the period in time for a
quasi-periodic flow when the observer sees it as periodic. (b) The appropriate value
of $c$ that converts a quasi-periodic flow in periodic, together with the
corresponding curve for periodic flows in lighter colors (cf. \reffig{c-per_1_1}).
The `*' correspond to the Hopf bifurcations at $Re_{p1}$ presented in
\reffig{bifur-per_1_1}}

In order to trace the curve $H_q(Re,c,U) = 0$, we utilize the same continuation
method as for system $H_p(Re,c,U)=0$ in \eqref{NSredE}, differing essentially in the
definition of the equation to vanish: for periodic flows the computations are much
simpler and faster than for quasi-periodic ones. The solution of \eqref{ceros} needs
an initial guess, which is obtained as described in the following subsections. Once
we have a quasi-periodic flow such that $H_q(Re,c,U) = 0$, we measure its amplitude
$A$, as in the case of periodic flows (cf. \eqref{amplitud}). If $u(x,y,t)$ is a
modulated wave, using \eqref{OM} and the $L$-periodicity we have
\begin{equation}\int_0^L \int_{-1}^1 u(x,y,\tau)^2 \dint y\dint x =
\int_0^L \int_{-1}^1 u(x-\phi,y,0)^2 \dint y\dint x = \int_0^L \int_{-1}^1
u(x,y,0)^2 \dint y\dint x.\label{Ampl_per}\end{equation} Since we have numerically
checked that $A$ is not constant for modulated waves, we conclude from
\eqref{Ampl_per} that it is a $\tau$-periodic function of $t$. This is not so for
the rotating waves of \refsecc{solper}, for which $A$ is constant in time. In the
case of a quasi-periodic flow $U(t)$, with the purpose of considering a concrete
value for the amplitude, we evaluate $A(t)$, at $t$ such that $U(t) \in \Sigma_1
\cap \Sigma_2$. This is simply a representative and easy to compute value for
$A(t)$, as we cannot obtain a single value for the amplitude of this class of flows.
We use it to trace the continuation curve: it provides the distance to the laminar
solution at some time instant. In the same way, $c$ can be considered as a
representative and time independent value for every quasi-periodic flow, so that we
can as well use it to trace continuation curves.

\subsection{Hopf bifurcation at $Re_{p1}$}\label{bif_Re_p1}

For $Re_p < Re_{p1}$ and $\alpha\approx1$ the corresponding time-periodic flow in
\reffig{bifur-perp_22x70} is unstable, but its temporal evolution ultimately decays
to the laminar flow. The results in \citet{Soibelman91} point out the existence of a
subcritical Hopf bifurcation at $Re_{p1}$: they use the vorticity equation and only
consider $N\leqslant2$ Fourier modes. According to bifurcation theory (see
\citet{Marsden76}), this means that the bifurcating quasi-periodic flows are locally
stable. Following this result we tried to find quasi-periodic flows in the
subcritical region, but with no success: we were not able to detect a quasi-periodic
attracting solution for $Re_p < Re_{p1}$ and $N\geqslant3$. In consequence we direct
the search of quasi-periodic flows to the supercritical region, i.e. for $Re_p >
Re_{p1}$.

It is also known from bifurcation theory that, in the case of a supercritical Hopf
bifurcation in which fixed points before it are unstable and after it stable, the
branch of periodic solutions that emanates at the bifurcation point is locally
unstable. That should be the situation for $Re_p > Re_{p1}$ and thus the bifurcating
quasi-periodic flows are locally unstable and therefore hard to locate by direct
numerical integration, because the evolution of the fluid close to them does not
remains near as time evolves.

On the other hand, pseudo-Newton's method applied to solve \eqref{ceros} does not
distinguish between stable or unstable flows. However, the difficult task is the
search of a good starting guess for \eqref{ceros}; it is carried out as follows.
First we consider an unstable time-periodic solution $U^p_1(t)$, i.e. a fixed point
of \eqref{NSred}, for certain $Re_1 < Re_{p1}$. Close to the Hopf bifurcation
$Re_{p1}$ in the unstable region $Re_p<Re_{p1}$, the linearization of \eqref{NSred}
around $U^p_1$ has just a couple of complex conjugate eigenvalues,
$\lambda=\lambda_r\pm i\lambda_i$ with $\lambda_r>0$ (cf. \refsecc{solper}). If
$w=w_r\pm iw_i$ is the associated eigenvector, we choose $v\in \langle
w_r,w_i\rangle$, i.e. $v$ is in the plane of the most unstable directions, along
which the flow escapes in the fastest fashion from $U^p_1$. Now we change to $Re_2 >
Re_{p1}$ close to the bifurcation and select $c_2$ near to the one associated with
the corresponding periodic orbit $U^p_2$ at $Re_2$, i.e. $H_p(Re_2,c_2,U^p_2)=0$.
For $Re_2, c_2$ and $|r|\ll1$ a small constant, we integrate numerically the initial
condition $U=U^p_1 + rv$, until a time when it become as closest as possible to a
solution of $H_q(Re_2,c_2,U)=0$. This first approximation of $c_2$ is optimized by
using minimization algorithms. In this case we also observe a value $\tau$ of the
return time of $P_c$, neighbouring to $2\pi/|\Imag(\lambda)|$, taking $\lambda$ as
the unique purely imaginary eigenvalue (together with its conjugate) at the Hopf
bifurcation $Re_{p1}$.

\figura{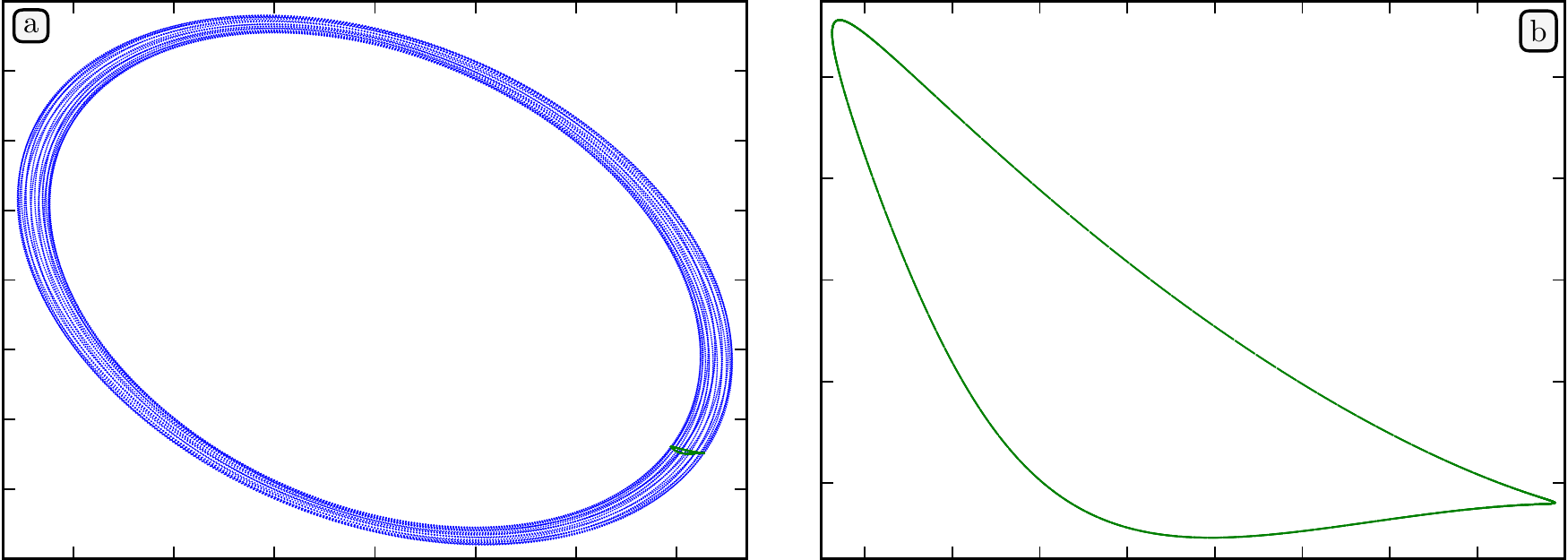}{Two representations of the solution vector $U(t)$ of an
(almost resonant) quasi-periodic flow projected on the plane of the same coordinates
as \reffig{secc-toro1}, for $Re_p=7000, \alpha=1.02056, N=22, M=70$. The range of
values for (a) is $[-0.00185,0.00185] \times [-0.004,0.004]$ and $[0.00147,0.00164]
\times [-0.002495,-0.002385]$ for (b). In (a) a dot is plotted each $\Delta t=0.01$
time units. In (b) a dot is plotted only for $t$ such that $U(t) \in \Sigma_1$. The
curve in (b) is also represented as a reference on the lower right corner of (a) as
a small green line}

The initial guess obtained in this way is improved through pseudo-Newton's
iterations applied to the function $H_q$, to finally obtain a first unstable
modulated wave for $Re_p > Re_{p1}$. In \reffig{secc-toro1} we present numerical
evidence that solutions $U(t)$ of \eqref{ceros} lie in a $2$-torus and are unstable
in time for supercritical $Re_p$. In a) we observe how the trajectory, projected on
a plane of two arbitrary coordinates of $U(t)$, fills densely a $2$-torus. On the
outer red curve of b), it is simply plotted $U(t)$ when $U(t) \in \Sigma_1$ with the
following adaptation: if $t > \tau$, being $\tau$ the period of $U(t)$ as a
modulated wave, we plot $U(s)$ for $s$ such that $0\leqslant s=t-p\tau<\tau$ and
$p=[t/\tau]$, ($[\cdot]$ stands for the integer part) i.e. we treat $U(t)$ as if it
were exactly $\tau$-periodic in order to avoid its unstability. We note that, as we
are on the Poincar\'e section $\Sigma_1$, the $2$-torus in a) is reduced to a closed
curve in b), which corresponds to the unstable quasi-periodic flow, seen as if it
were unperturbed by numerical errors. Likewise in b) we let the flow evolve for long
time and plot again $U(t)$ (in green) when $U(t) \in \Sigma_1$ but suppressing the
previous time adaptation. Here we can check that the flow is unstable, because it
moves away from the outer closed curve and falls to the time-periodic and stable
r\'egime: a simple point in the centre of the figure.

\figura{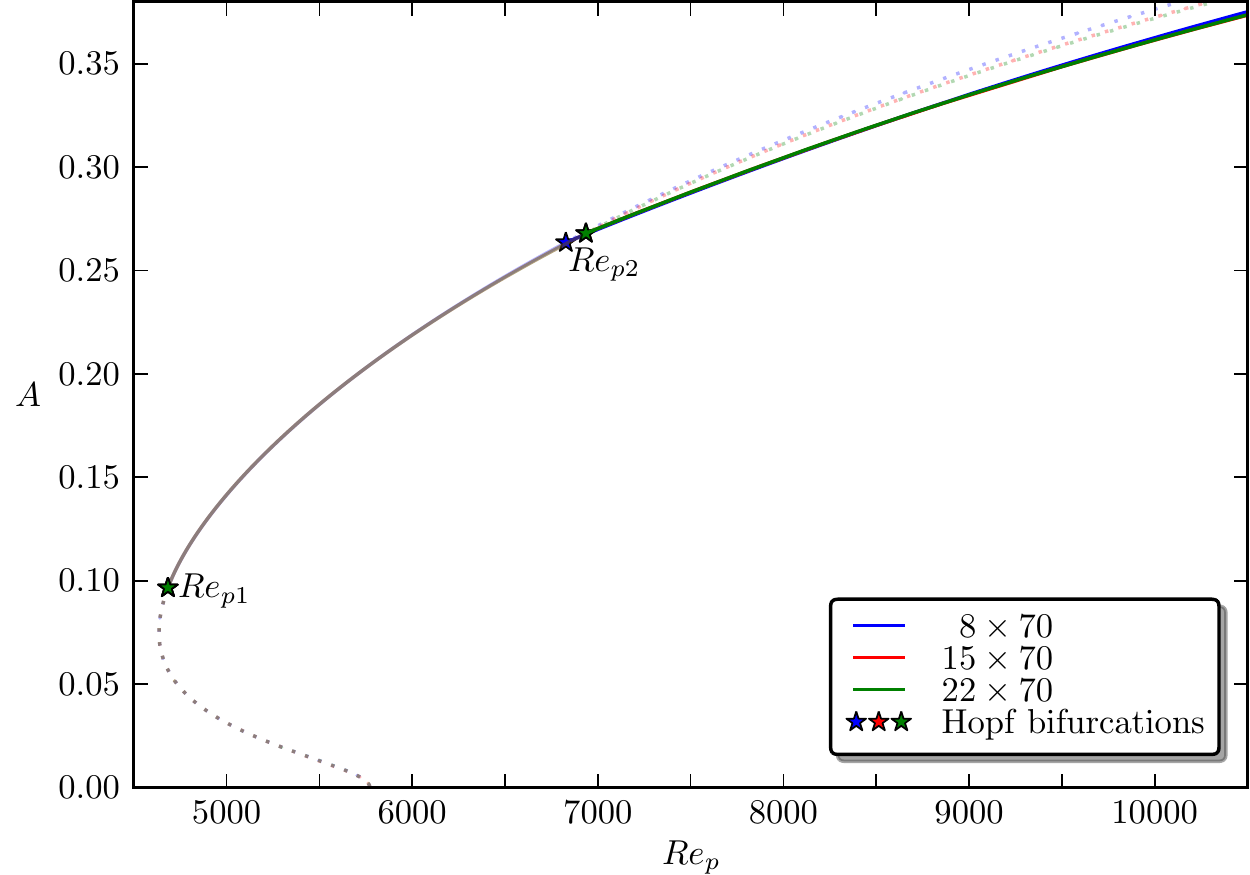}{Bifurcated branches of quasi-periodic flows at the
supercritical Hopf bifurcation of periodic flows at $Re_{p2}$. Each curve represents
$A$ as a function of $Re_p$ and has specified $N\times M$. Both branches of periodic
(in lighter colors) and quasi-periodic flows are presented (the green color almost
completely conceals the red one). Calculations are shown for $\alpha=1.02056$ and
$\Delta t=0.02$. The `*' correspond to the $Re_{p1},Re_{p2}$ Hopf bifurcations
presented in \reffig{bifur-perp_22x70}. For the range shown, the bifurcating branch
consist of quasi-periodic orbits, stable to disturbances of the same wavenumber}

\figura{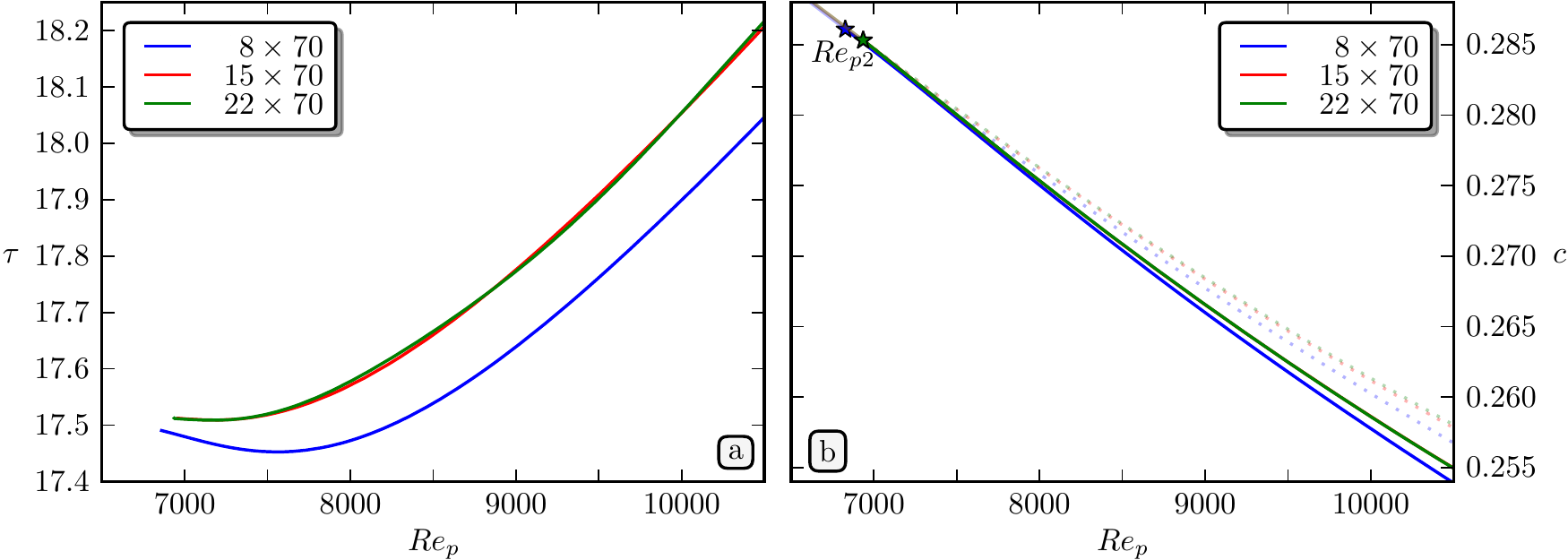}{Different curves related to \reffig{bifur-cper-Rep2_1_02}
using $Re_p$ in the abscissa axis for $N\times M$ points as specified. (a) and (b)
are as in \reffig{c-tau-Rep1_1_1}. In (b), in lighter colors, it is also drawn the
periodic flows around $Re_{p2}$ (marked with `*'). In this case, the green color
almost completely conceals the red one}

Once we have obtained the first solution of \eqref{ceros} by the pseudo-Newton's
method, we use continuation methods to traverse the bifurcating branch of
quasi-periodic flows parametrized by $Re_p$. In \reffig{bifur-cper1_1_1} we plot
the amplitude $A$ for each quasi-periodic solution as a function of $Re_p$. It seems
that we have achieved both qualitative and quantitative convergence because we
obtain a similar graph in increasing the values of $N,M$. It looks also clear from
this plot that the Hopf bifurcation is supercritical so, at least locally, solutions
on the bifurcating branch are unstable. The analysis of the stability of a
quasi-periodic solution is done by means of the eigenvalues of the linear part of
$P_c$ at a fixed point $U(t)$ such that $P_c(U(0))=U(0)$. This computation,
analogously to the matrix $DH_q$, is obtained by extrapolated finite differences. In
the range of $Re_p$ presented in \reffig{bifur-cper1_1_1}, the quasi-periodic
solutions are unstable. Furthermore, unlike the bifurcation at $Re_{p2}$, solutions
at $Re_{p1}$ present certain symmetry: $\hat u_k(y,t)$ (defined in
\refsecc{numerical_app}) is an even or odd function of $y$ according to the parity
of $k$.

We also observe in \reffig{c-tau-Rep1_1_1}(a) for different values of $N$, an
indicator of the numerical effort involved in the evaluation of the map $P_c$. The
big slope of the $Re_p$-$\tau$ curve shows the high computational cost involved as
$Re_p$ is slightly increased. For $Re_p\in [3840,3872]$, $\alpha = 1.1$ and $\Delta
t=0.02$ the time needed to return to $\Sigma_1$ is $\tau \in [7000,10700]$. Likewise
it has been necessary to adapt dynamically the minimum number of times, $n_c$
(defined just before system \eqref{ceros}), that the solution crosses $\Sigma_1$
before it returns to the starting point. For $Re_p$ close to $Re_{p1}$, $n_c$ starts
at $1$ and it is incremented by $1$ when $Re_p$ varies approximately in just
$2$--$3$ units. The way we modify $n_c$ is described next. If we change slightly
$Re_p$ the return time $\tau$ should be varied in accordance. Thus, we impose that
the new value obtained for $\tau$ satisfy $|\tau - \tau_o| < \varepsilon$, for some
tolerance $\varepsilon$ and $\tau_o$ the former return time. If this condition is
not fulfilled we increase or decrease $n_c$ by one unit, until it is satisfied, or
we decrease $Re_p$ if necessary.

In \reffig{c-tau-Rep1_1_1}(b) is represented $c$ for the same range of $Re_p$. We
observe for $c$, nearby values as their counterpart periodic flows (cf.
\reffig{c-per_1_1}), but in this case the curve has a large decreasing slope. We
emphasize that the $c$ graph shows a bifurcation diagram of periodic and
quasi-periodic flows, which is time independent in contrast to the $Re$-$A$ plot in
\reffig{bifur-cper1_1_1}.

\subsection{Bifurcation at $Re_{p2}$}\label{bif_Re_p2}

\figura{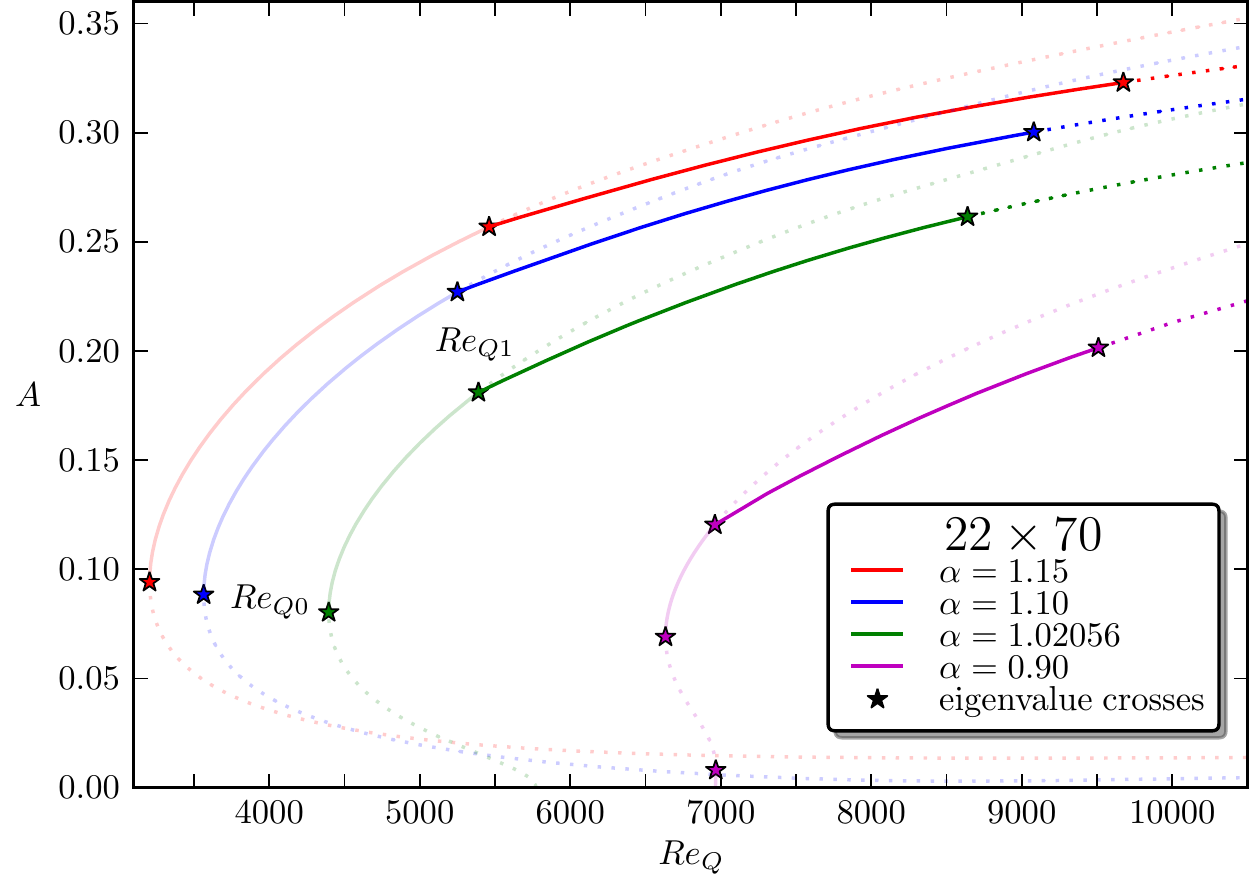}{Analogous to \reffig{bifur-cper-Rep2_1_02} based on
$Re_{Q1}$, using the specified $\alpha$ and $N\times M = 22 \times70$, $\Delta t =
0.01$. The bifurcated branches of quasi-periodic solutions have a change of
stability at another Hopf bifurcation. Continuous and discontinuous lines represent
stable and unstable flows respectively}

\figura{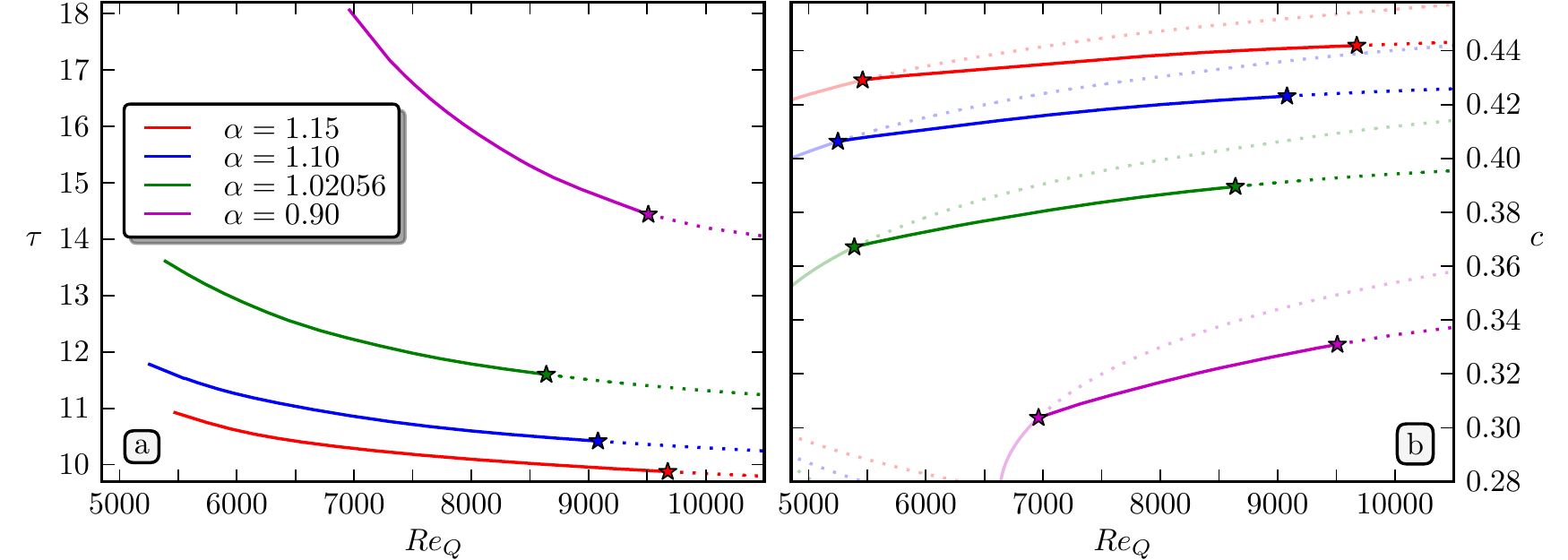}{Values of $\tau$ (a) and $c$ (b) associated to quasi-periodic
flows of \reffig{bifur-cper-ReQ1}. In (b), in lighter colors, it is also drawn the
periodic flows around $Re_{Q1}$ of \reffig{bifur-cper-ReQ1}. Colors for curves in
(a) are also valid for (b)}

\figura{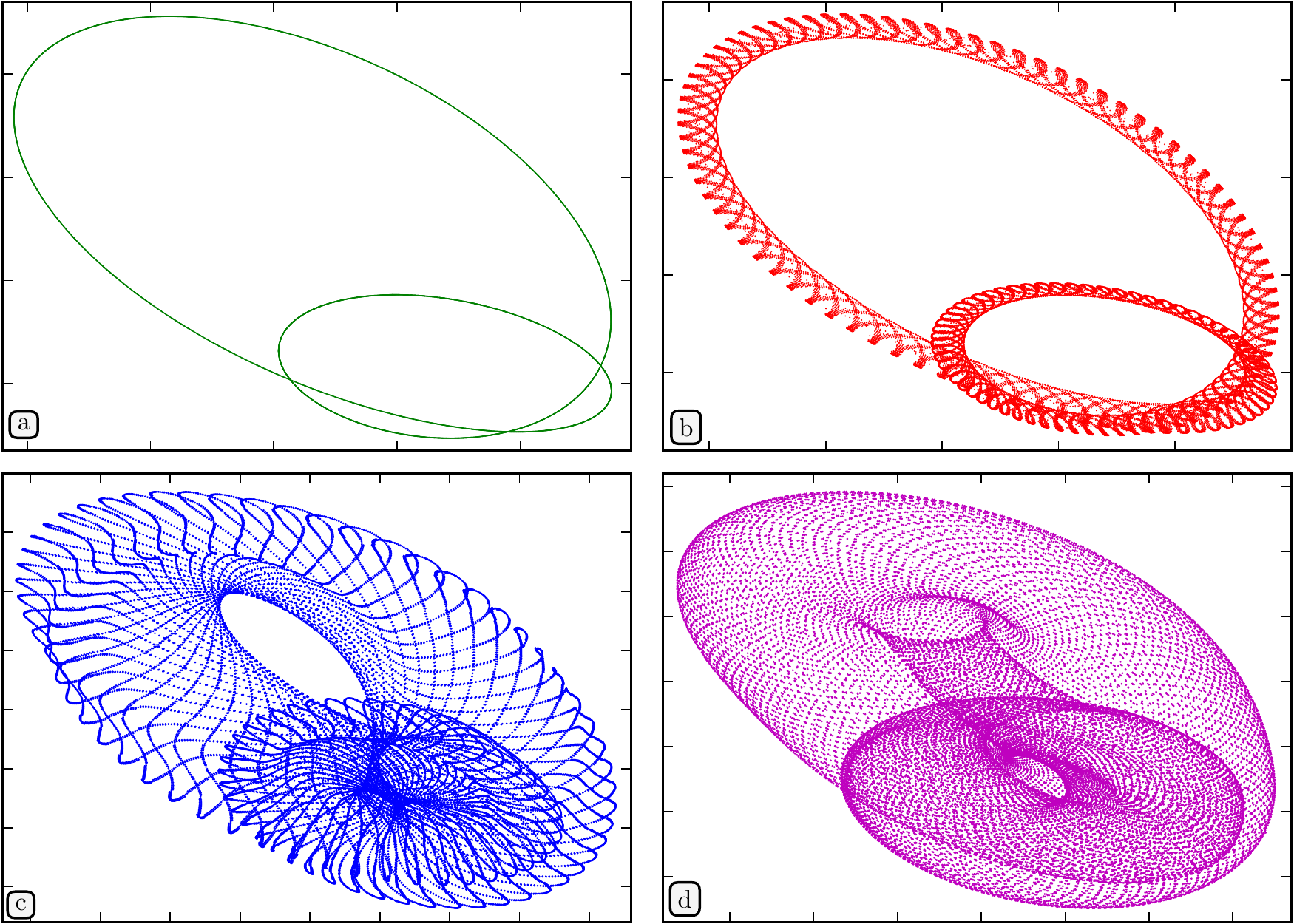}{(a) An stable $2$-torus on $\Sigma_1$ for $Re_Q = 8635$,
$\alpha = 1.02056$, $N\times M = 22 \times70$, $\Delta t = 0.009$. The integration
time on the figure is about $1{,}222{,}000$ time units. (b) Analog of (a) for $Re_Q
= 8640$. We observe in this case that the initial $2$-torus is unstable and is
attracted by a $3$-torus. The integration time on the figure is about $419{,}000$
time units. (c) The same as (b), but the initial condition is an unstable $2$-torus
for $Re_Q=9750$ which is also attracted by a $3$-torus. The integration time is
about $1{,}005{,}000$ time units. (d) The initial condition is an unstable $2$-torus
for $Re_Q=10500$ and $\alpha=1.10$ which is again attracted by a $3$-torus. The
integration time is about $1{,}002{,}000$ time units. The respective ranges in the
four figures are $[0.0019, 0.00445] \times [-0.00965, -0.0053]$, $[0.0018, 0.0045]
\times [-0.0098, -0.0052]$, $[0.0013, 0.0058] \times [-0.0126, -0.005]$ and
$[0.0016, 0.00535] \times [-0.0147, -0.0078]$}

\begin{figure}[ht]\begin{center}
\includegraphics{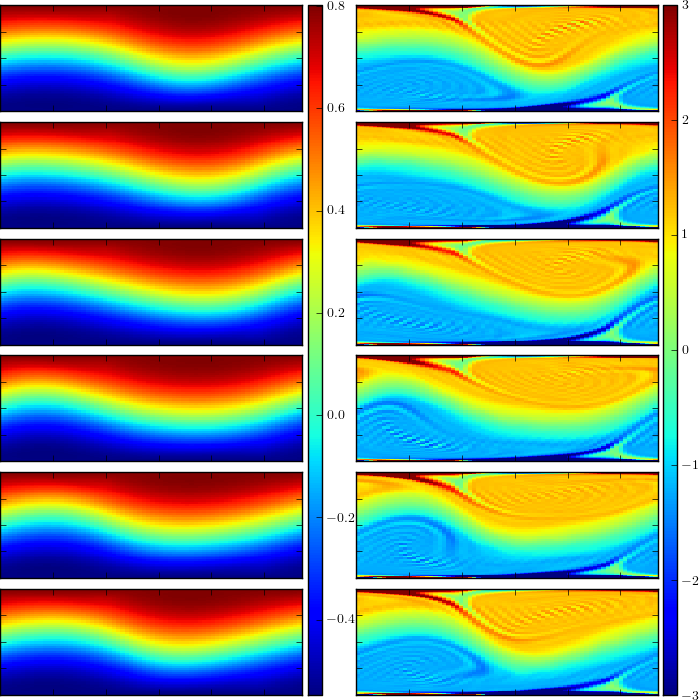}

\caption{Stream lines (left) and levels of vorticity (right) for a unstable
quasi-periodic flow at $Re_Q = 10500, \alpha=1.1, N=22, M=70, \Delta t=0.009$. Surface
levels are plotted for $0\leqslant x\leqslant L, -1\leqslant y\leqslant1$. The
depicted time instants from top to down are $0, \tau/6, 2\tau/6, 3\tau/6, 4\tau/6$,
and $5\tau/6$ for the associated orbit period $\tau=10.24$. Ranges for levels of
stream lines and vorticity are $[-0.55, 0.80]$ and $[-3, 3]$ respectively.
They are represented by a common scale of colors to the right of each column.}
\label{peli-vorti}\end{center}\end{figure}

As we presented in \reffig{bifur-perp_22x70}, for $Re > Re_{p2}$ the corresponding
periodic flow is unstable, so the evolution of \eqref{NSred} from such a flow as
initial condition, drives the fluid away from it. By following the temporal
evolution of this flow, we observe that the fluid seems to fall in a regular
r\'egime, which finally proves to be a quasi-periodic attracting solution. This is
checked in \reffig{secc-toro-Rep2} where we plot the projection of the solution
vector $U(t)$ over the plane of the same two coordinates as in \reffig{secc-toro1}.
Each point in \reffig{secc-toro-Rep2}(a) corresponds to the value of the specified
coordinates at a time instant. As we can observe, the trajectory appears to fill
densely the projection of a $2$-torus. Plotting the same coordinates as above, but
only when $U(t) \in \Sigma_1$, we see in \reffig{secc-toro-Rep2}(b) a closed curve,
which seems again to confirm that the flow lives in a $2$-torus.

Let us denote as $U^0(t) = (\bar u_0^0, \dots,\bar u_N^0)(t)$, a time instant of the
flow in the attracting $2$-torus. We need to approximate the value of $c^0$ which
better makes $U^0$ appear as a periodic flow. That $c^0$ exists according to
\eqref{OM}--\eqref{RedaPer}. It can be estimated as (cf. \citet{Rand82})
\[c^0 = \frac{2\pi}{k}\lim_{t\rightarrow\infty}
\frac{n(t)}{t},\quad\mbox{for}\quad n(t)=\left[\frac{\arg(\Conj{(\bar u^0_{km}(t))})
  }{2\pi} \right],\] where $k$ is the number of peaks of the wave $U^0$ for
$x\in[0,L]$ and we have used the midpoint of the channel for $m=M/2$. With $t$ large
enough, $2\pi n(t) / kt$ gives an approximation of $c^0$ which we optimize using
minimization. Next we use $U^0, c^0$ as the initial condition for pseudo-Newton's
method applied to \eqref{ceros} to confirm that our attracting 2-torus $U^0$ is in
fact a modulated wave. For a fixed $\alpha$, once we have a first point $(Re_p, c^0,
U^0)$ which satisfies \eqref{ceros}, taking $Re_p$ as a continuation parameter, we
can trace the curve of quasi-periodic flows in the $Re$--$A$ plane by
pseudo-arclength numerical continuation, applied to $H_q$ as in \refsecc{bif_Re_p1}.

As we observe in \reffig{bifur-cper-Rep2_1_02}, there appears a branch of
quasi-periodic solutions which bifurcates supercritically from the curve of periodic
flows. We may again have zeros of $H_q(Re,c,U)$ that can either correspond to stable
or unstable quasi-periodic solutions. Since on crossing the bifurcation point
$Re_{p2}$, the periodic orbits change from stable to unstable, the branch of
bifurcating quasi-periodic flows are locally stable to two-dimensional superharmonic
disturbances. By means of the eigenvalues of the Jacobian matrix $\PA P_c/\PA u$ we
also compute the stability of the obtained quasi-periodic solutions. For
$\alpha=1.02056$ and the range of $Re_p\in[7000,13000]$ studied, all quasi-periodic
flows found are stable to perturbations of the same wavelength and the situation is
kept when $N,M$ are increased. In \reffig{c-tau-Rep2_1_02} we present the curves of
frequencies $c, \tau$ which define the different modulated waves. Again the $Re$-$c$
graph shows a time independent bifurcation diagram.

\subsection{Bifurcation at $Re_{Q1}$}\label{bif_Re_Q1}

In the case of constant flux the bifurcation diagram of periodic flows is
qualitatively different to that of constant pressure, as can be verified in
\reffig{bifur-per_1_1}. For $Re_Q$ and the values of $\alpha$ considered ($0.9$,
$1.02056$, $1.1$ and $1.15$) there is a change of stability at the minimum Reynolds
$Re_{Q0}$ of the amplitude curves, but no new bifurcations are born there. The first
Hopf bifurcation occurs at the point labeled $Re_{Q1}$ in \reffig{bifur-per_1_1}. We
can summarize that the qualitative picture of amplitudes of quasi-periodic solutions
emanating from $Re_{Q1}$ is analogous to that of $Re_{p2}$. The main differences are
basically quantitative, because $Re_{Q1} < Re_{p2}$ (cf. figures
\ref{bifur-cper-Rep2_1_02} and \ref{bifur-cper-ReQ1}).

For $\alpha = 0.9$, $1.02056$, $1.1$ and $1.15$, we compute the quasi-periodic flows
that bifurcate from $Re_{Q1}$, following the same steps of \refsecc{bif_Re_p2}. In
\reffig{bifur-cper-ReQ1} we plot their amplitudes, and again, as in the case of
$Re_{p2}$ we observe a supercritical bifurcation. The associated frequencies $c$,
$\tau$, to the modulated waves are presented in \reffig{c-tau-ReQ1}. We remark the
analogies between bifurcation diagrams of $A$ and $c$ in respective figures
\ref{bifur-cper-ReQ1} and \ref{c-tau-ReQ1}(b), the latter one being time
independent. The quasi-periodic solutions found from $Re_{Q1}$ and $\alpha=1.02056$,
are stable for $Re_{Q1} < Re_Q \lesssim 8640$. At $Re_Q \approx 8640$ the branch of
quasi-periodic solutions loses stability at a new Hopf bifurcation, giving rise to a
family of attracting tori of $3$ frequencies. Numerical evidence of this bifurcation
is shown in \reffig{secc-toro3d-ReQ}(a) where the same coordinates of $U(t)$ as in
\reffig{secc-toro1} are plotted on $\Sigma_1$ for $Re_Q = 8635$, yielding an
apparently perfect closed curve, after a considerably long time of integration: it is
a stable quasi-periodic solution. On the contrary, in a similar plot for $Re_Q =
8640$, \reffig{secc-toro3d-ReQ}(b) shows an unstable quasi-periodic flow, which is
attracted by a $3$-torus. The new frequency of this attracting solution is verified
using in turn the Poincar\'e section $\Sigma_2$. We have carried this out by
plotting the same two selected coordinates of $U(t)$ only when the flow crosses
$\Sigma_1$ if, in addition, it is approximately on $\Sigma_2$. We obtain in this way
what seems to be a closed curve (see \citet{Casas04} for a plot) as may be expected
for a $3$-torus. We can observe two other more involved $3$-torus in
\reffig{secc-toro3d-ReQ}(c) and (d).

Finally in \reffig{peli-vorti}, stream lines and levels of vorticity of an unstable
quasi-periodic flow for $Re_Q=10500$, $\alpha=1.1$, $N\times M=22\times 70$ and
$\Delta t=0.009$, are presented at six equidistant values of time in the interval
$[0,\tau]$ for $\tau=10.24$. We observe (as expected) larger vorticity close to the
walls than in the channel centre, which in addition can be confirmed on the stream
lines figures.

\section{Conclusions}\label{conclusions}

In this work we have studied some bifurcations of plane Poiseuille flow. We have
reproduced results of other authors and obtained similar qualitative results about
the Hopf bifurcations, in what concerns to their number and location. The main
quantitative differences between \citeauthor{Soibelman91}'s \cite{Soibelman91}
computations and ours are due to the larger resolution we have used, together with
the distinct formulations implemented of the Navier--Stokes equations. The important
qualitative difference is the kind of bifurcation found at $Re_{p1}$: in their
computations this bifurcation is subcritical, but improving the precision of the
numerical approach we obtain that it is supercritical. Then, the bifurcating
quasi-periodic orbits are unstable. This has also been confirmed by numerical
simulations.

In the case of the bifurcation at $Re_{p1}$, because the lengthy time integrations,
we have only been able to move away a few tens from $Re_{p1}$. The further we
advanced in $Re_p$, the greater are the numerical difficulties we encounter to track
the bifurcating branch of quasi-periodic flows, due to long time integrations. It is
also worth to mention the complications derived from their instability. Close to
$Re_{p1}$ and with the discretization employed ($N=8,\ M=70,\ \Delta t=0.02$), it
seems that we have achieved both qualitative and quantitative convergence. By
observing \reffig{bifur-cper1_1_1} we can conjecture that, for the range of
$\alpha\in [1,1.1]$ considered, the minimum $Re\approx 2900$ attained with
travelling waves is not lowered by quasi-periodic flows. This question still remains
open for two-dimensional flows, although \citet{Ehrenstein91} solved the gap between
experiments and numerical results in the case of three-dimensional flows. However,
in the present work for $\alpha=0.89$, $Re_{p1}\approx 7250$ we have localized a
subcritical branch of stable quasi-periodic orbits. They are difficult to follow by
means of the approach described in \refsecc{bif_Re_p1}, because the time needed to
evaluate the Poincar\'e map is $\tau\approx 15,000$ time units. It remains open
whether that family could reduce the minimum $Re\approx 2900$ of periodic flows to
$Re\approx 1000$, where transition has been observed experimentally.

For $Re_p>Re_{p2}$ the quasi-periodic flows encountered are attracting and the
integration time is of the order of tens, so in this case the computational cost is
drastically reduced compared with the bifurcation at $Re_{p1}$. The range of $Re_p$
obtained for attracting quasi-periodic flows moves now to several thousands.
However, in spite of keeping qualitative convergence, the use of larger Reynolds
numbers makes necessary an increase in precision to get, furthermore, quantitative
convergence. An analogous qualitative picture is found at the quasi-periodic flows
which bifurcate from $Re_{Q1}$. In this case the quasi-periodic solutions quickly
loses stability and we have also obtained another Hopf bifurcation to a family of
tori with $3$ basic frequencies. We could say that dynamics are richer for $Re_Q$
than $Re_p$ (see \citet{Casas04}), because bifurcations and different vortical
states appear for lower $Re_Q$ than the counterpart $Re_p$.

As future work, it would be of interest to analyse the stability to disturbances
with different wavenumber $\alpha$ or even to $3$-dimensional perturbations. The
connections of the different families of solutions is also of great relevance, or
even the discovering of new vortical states which could approach more the transition
to turbulence. Likewise, due to nonnormality in the Navier-Stokes system, the
sensitivity of eigenvalues to perturbations is an important issue that can be
analyzed by means of the pseudospectra. This study would give a measure of the
reliability of the spectrum obtained in linearizing \eqref{NSred} around periodic
flows, mainly for high values of $Re$. In this work, the stability according to
eigenvalues is coherent with our direct numerical simulations \eqref{IN} of the
different flows.

\begin{acknowledgments}
We thank C. Sim\'o and J. Sol\`a-Morales for valuable discussions during the
preparation of this paper. P.S.C. has been partially supported by funds from the
Departamento de Ma\-te\-m\'a\-ti\-ca Aplicada I (Universidad Polit\'ecnica de
Catalu\~na), and the MCyT-FEDER grant MTM2006-00478. A.J. has been supported by
the MEC grant MTM2009-09723 and the CIRIT grant 2009SGR-67. The computing facilities
of the UB-UPC Dynamical Systems Group (clusters Hidra and Eixam) have been widely
used.
\end{acknowledgments}

\bibliography{poiseuille}
\nocite{Drissi99,Herbert76,Herbert91,Pugh88,Rozhdestvensky84,Zahn74}

\end{document}